\documentclass[reqno]{amsart}

\usepackage{amsmath}
\usepackage{amsfonts,amssymb}
\usepackage[dvips]{graphicx}
\usepackage{epsfig}
\newtheorem{theorem}{Theorem}[section]
\newtheorem{remark}{Remark}[section]
\newtheorem{lemma}[theorem]{Lemma}
\newtheorem{definition}{Definition}[section]

\newtheorem{proposition}[theorem]{Proposition}
\numberwithin{equation}{section}
\begin{document}

\title[the compressible Euler equations with a time periodic outer force]
	{Existence of a time periodic solution for the compressible Euler equation with a time periodic outer force}
%\subtitle{}
\author{Naoki Tsuge}
\address{Department of Mathematics Education, 
Faculty of Education, Gifu University, 1-1 Yanagido, Gifu
Gifu 501-1193 Japan.}
\email{tuge@gifu-u.ac.jp}
\thanks{
N. Tsuge's research is partially supported by Grant-in-Aid for Scientific 
Research (C) 17K05315, Japan.
}
\keywords{The Compressible Euler Equation, a time periodic outer force, the compensated compactness, a time periodic solution, the modified Lax Friedrichs scheme, the fixed point theorem.}
\subjclass{Primary 
35L03, 
%Initial value problems for first-order hyperbolic equations,
35L65, 
%Conservation laws,
35Q31, 
%Euler equations,
76N10,
%Existence, uniqueness, and regularity theory, 
76N15; 
%Gas dynamics, general; 
Secondary
35A01, 
%Existence problems: global existence, local existence, non-existence, 
35B35,   
%Stability 
35B50, 
%Maximum principles,
35L60,   
%Nonlinear first-order hyperbolic equations, 
76H05,   
%Transonic flows,
76M20.   
%Finite difference methods. 
}
\date{}

\maketitle

\begin{abstract}
We are concerned with a time periodic supersonic flow through a bounded interval.
This motion is described by the compressible Euler equation with a time periodic 
outer force. Our goal in this paper is to prove the existence of a time periodic 
solution. Although this is a fundamental problem for other equations, it has not been 
received much attention for the system of conservation laws until now.

When we prove the existence of the time periodic solution, we face with two problems. 
One is to prove that initial data and the corresponding solutions after one period 
are contained in the same bounded set. To overcome this, we employ the generalized invariant region, which depends on the space variables. This enable us to investigate the 
behavior of solutions in detail. Second is to construct a continuous map.
We apply a fixed point theorem to the map from
initial data to solutions after one period. 
Then, the map needs to be continuous. To construct this, we introduce 
the modified Lax-Friedrichs scheme, which has a recurrence formula consisting of 
discretized approximate solutions. The formula yields 
the desired map. Moreover, the invariant region grantees that
it maps a compact convex set to itself. In virtue of the fixed point theorem, we can prove a existence 
of a fixed point, which represents a time periodic solution.
Finally, we apply the compensated compactness framework to prove the convergence of our approximate solutions. 
\end{abstract}

%\tableofcontents

\section{Introduction}
The present paper is concerned with the compressible Euler equation 
with an outer force. 
\begin{align}
\begin{cases}
\displaystyle{\rho_t+m_x=0,}\\
\displaystyle{m_t+\left(\frac{m^2}{\rho}+p(\rho)\right)_x
	=F(x,t)\rho,}
\end{cases}x\in(0,1),\quad t\in(0,1)
\label{eqn:force}
\end{align}
where $\rho$, $m$ and $p$ are the density, the momentum and the 
pressure of the gas, respectively. If $\rho>0$, 
$v=m/\rho$ represents the velocity of the gas. For a barotropic gas, 
$p(\rho)=\rho^\gamma/\gamma$, where $\gamma\in(1,5/3]$ is the 
adiabatic exponent for usual gases. The given function 
$F\in C^1([0,1]\times[0,1])$ represents a time periodic outer force with the time period $1$, 
i.e., $F(x,0)=F(x,1)$.

We consider the initial boundary value problem (\ref{eqn:force}) 
with the initial and boundary data
\begin{align}  
(\rho,m)|_{t=0}=(\rho_0(x),m_0(x))\quad (\rho,m)|_{x=0}=(\rho_b,m_b).
\label{eqn:I.C.}
\end{align}
The above problem \eqref{eqn:force}--\eqref{eqn:I.C.} can be written in the following form 
\begin{align}\left\{\begin{array}{lll}
u_t+f(u)_x=g(x,t,u),\quad{x}\in(0,1),\quad t\in(0,1),\\
u|_{t=0}=u_0(x),\\
u|_{x=0}=u_b
\label{eqn:IP}
\end{array}\right.
\end{align}
by using  $u={}^t(\rho,m)$, $\displaystyle f(u)={}^t\!\left(m, \frac{m^2}{\rho}+p(\rho)\right)$ and 
$\displaystyle{g(x,u)={}^t\!\left(0,F(x,t)\rho\right)}$. We shall consider solutions with positive characteristic 
speeds. Therefore, from the Lopantinski condition, we do not supply a boundary condition at $x=1$.

In the present paper, we consider the compressible Euler equation. Let us survey the related mathematical results.

Concerning the one-dimensional initial value problem, {\sc DiPerna} \cite{D1}
proved the global existence by the vanishing viscosity method and a compensated compactness argument. 
The method of compensated compactness was introduced
by {\sc Murat} \cite{M1} and {\sc Tartar} \cite{Ta1,Ta2}.
{\sc DiPerna} first applied the method to systems for 
the special case where $\gamma=1+2/n$ and $n$ is an odd
integer. Subsequently, {\sc Ding}, {\sc Chen} and {\sc Luo} \cite{DC1} 
and {\sc Chen} \cite{C2} extended his analysis to any $\gamma$ in $(1,5/3]$. In \cite{DC2}, {\sc Ding}, {\sc Chen} and {\sc Luo} treated isentropic gas dynamics with a source term by using the fractional step procedure. By this result, the global existence theorem was obtained for the compressible Euler equation with an outer force. On the other hand, the stability (i.e., solutions are contained in a bounded set independent of the time variable) has not yet proved for a long time. Recently, in \cite{T6} and \cite{T9}, it is proved by using the generalized invariant region. This method is also employed in \cite{T2}--\cite{T8}.

In this paper, we consider a time periodic outer force and prove the existence of a time periodic solution.
Although this problem is fundamental for other equations (ex. \cite{MN}), it has not been received much attention until now. For a single conservation law, Takeno \cite{Takeno} proved the existence of a time periodic solution for the space periodic boundary condition. However, unfortunately, a little is known for the system of conservation laws.
This is reason why there are the following two problems. 
One is to prove that initial data and the corresponding solutions after one period 
are contained in the same bounded set. To overcome this, we employ the generalized invariant region, which depends on the space variables. Second is to construct a continuous map in a finite dimension.
We apply a fixed point theorem to a map from
initial data to solutions after one period. 
Then, the map needs to be continuous. To construct this, we introduce 
the modified Lax-Friedrichs scheme, which has a recurrence formula \eqref{eqn:recurrence1} consisting of 
discretized approximate solutions. It yields 
the continuous map. In addition, the invariant region grantees that
it maps a compact convex set to itself. Applying the Brouwer fixed point theorem, we can prove a existence 
of a fixed point, which represents a time periodic solution.

To state our main theorem, we define the Riemann invariants $w,z$, which play important roles
in this paper, as
\begin{definition}
	\begin{align*}
	w:=\frac{m}{\rho}+\frac{\rho^{\theta}}{\theta}=v+\frac{\rho^{\theta}}{\theta},
	\quad{z}:=\frac{m}{\rho}-\frac{\rho^{\theta}}{\theta}
	=v-\frac{\rho^{\theta}}{\theta}\quad
	\left(\theta=\frac{\gamma-1}{2}\right).
	%\label{eqn:Riemann-invariant}
	\end{align*}
\end{definition}
These Riemann invariants satisfy the following.
\begin{remark}\label{rem:Riemann-invariant}
	\normalfont
	\begin{align}
	&|w|\geqq|z|,\;w\geqq0,\;\mbox{\rm when}\;v\geqq0.\quad
	|w|\leqq|z|,\;z\leqq0,\;\mbox{\rm when}\;v\leqq0.\label{eqn:inequality-Riemann}\\
	&v=\frac{w+z}2,
	\;\rho=\left(\frac{\theta(w-z)}2\right)^{1/\theta},\;m=\rho v.
	\label{eqn:relation-Riemann}
	\end{align}From the above, the lower bound of $z$ and the upper bound of $w$ yield the bound of $\rho$ and $|v|$.
\end{remark}

Moreover, we define the entropy weak solution.
\begin{definition}
	A measurable function $u(x,t)$ is called an time periodic {\it entropy weak solution} of the initial boundary value problem \eqref{eqn:IP} with period $1$ if 
	\begin{align*}
&\int^{1}_{0}\int^{1}_0\rho(\varphi_1)_t+m(\varphi_1)_xdxdt+\int^{1}_{0}\rho_0(x)\left
(\varphi_1(x,0)-\varphi_1(x,1)\right)dx\\
&\quad+\int^{1}_{0}m_b\varphi_1(0,t)dt
=0,\\
&\int^{1}_{0}\int^{1}_0m(\varphi_2)_t+\left(\frac{m^2}{\rho}+p(\rho)\right)
(\varphi_2)_x+F(x,t)\rho\varphi_2 dxdt\\
&\quad+\int^{1}_{0}m_0(x)\left(\varphi_2(x,0)-\varphi_2(x,1)  \right)dx+\int^{1}_{0}\left\{\frac{\left(m_b\right)^2}{\rho_b}+p(\rho_b)\right\}\varphi_2(0,t)dt=0
\end{align*}
holds for any test function $\varphi_1,\varphi_2\in C^1_0([0,1]\times[0,1])$ satisfying 
$\varphi_1(1,t)=\varphi_2(1,t)=0$ and 
	\begin{align*}
	&\int^{1}_0\int^{1}_0\hspace{-1ex}\eta(u)\psi_t+q(u)\psi_x+\nabla\eta(u) g(x,u)\psi dxdt\geqq0
	\end{align*}
	holds for any non-negative test function $\psi\in C^1_0((0,1)\times(0,1))$, where 
	$(\eta,q)$ is a pair of convex entropy--entropy flux of \eqref{eqn:force}.
\end{definition}
\begin{remark}\normalfont
	It follows from the above definition that $f(u)$ converges in weak sense to $f(u_b)$, i.e., 
	\begin{align*}
	\frac1{\varepsilon}\int^{\varepsilon}_{0}f(u)\,dx\rightarrow f(u_b)\quad
	\text{    as $\varepsilon\rightarrow0$ in $L^{\infty}({\bf R}_+)$ weak${}^*$}.
	\end{align*}
\end{remark}

We choose a positive constant $K$ such that  
\begin{align}
|F(x,t)|\leqq K\quad x\in[0,1].
\label{eqn:condition-X}
\end{align}

We then choose a positive value $L$ and $M$ such that 
\begin{align}
M\geqq L\geqq\left(1+K\right).
\label{eqn:condition-M}
\end{align}

Then our main theorem is as follows.
\begin{theorem}\label{thm:main}
	
	If
	$u_0=({\rho}_0, {m}_0)\in{L}^{\infty}(0,1)$ satisfy
	\begin{align}
	%\begin{split}
	0\leqq\rho_0(x)
	,\;\; L-Kx\leqq{z}(u_0(x)),\;\; w(u_0(x))\leqq M+Kx,\quad x\in(0,1)
	%\end{split}
	\label{eqn:IC}
	\end{align}
	and the boundary data satisfy 
	\begin{align}
	L\leqq v_b-(\rho_b)^{\theta}/\theta,\quad v_b+(\rho_b)^{\theta}/\theta\leqq M,
	\label{eqn:BC}
	\end{align}
	%in terms of Riemann invariants, or
	%\begin{align*}
	%&0\leqq\rho_0(x),\quad
	%\int^{x}_0X(y)dy+\int^t_0T(s)ds-M
	%\leqq v_0(x)-\frac{\{\rho_0(x)\}^{\theta}}{\theta},\\
	%&v_0(x)+\frac{\{\rho_0(x)\}^{\theta}}{\theta}\leqq 
	%\int^{x}_0X(y)dy+\int^t_0T(s)ds+M
	%\end{align*}
	%in the physical variables. 
	then the initial boundary value problem $(\ref{eqn:IP})$ has a 
	time periodic entropy weak solution.

	Moreover, the solution satisfies 
	\begin{align}
	\begin{split}
	&0\leqq\rho(x,t)
	,\; L-Kx\leqq{z}(u(x,t)),\; w(u(x,t))\leqq M+Kx,\quad(x,t)\in(0,1)\times(0,1)
	\end{split}
	\label{eqn:bound}
	\end{align}
\end{theorem}

\begin{remark}
	\normalfont
	If solutions satisfy \eqref{eqn:bound}, both of their characteristic speeds are positive.
\end{remark}

\subsection{Outline of the proof}

The proof of main theorem is a little complicated. Therefore, 
before proceeding to the subject, let us grasp the point of the main estimate by a formal argument. 
We assume that a solution is smooth and the density is nonnegative in this section.

We consider the physical region $\rho\geqq0$ (i.e., $w\geqq z$.). Recalling Remark \ref{rem:Riemann-invariant}, it suffices to 
derive the lower bound of $z(u)$ and the upper bound of $w(u)$ to obtain the bound of $u$. To do this, we diagonalize \eqref{eqn:force}. 
If solutions are smooth, we deduce from \eqref{eqn:force} 
\begin{align}
z_t+\lambda_1z_x=F(x,t),\quad
w_t+\lambda_2w_x=F(x,t),
\label{eqn:force2}
\end{align} 
where $\lambda_1$ and $\lambda_2$ are the characteristic speeds defined as follows 
\begin{align}
\lambda_1=v-\rho^{\theta},\quad\lambda_2=v+\rho^{\theta}.
\label{eqn:char}
\end{align}

Moreover, set
\begin{align*} \tilde{z}=z+Kx,\;\tilde
{w}=w-Kx. 
\end{align*}
Then, it 
follows from \eqref{eqn:force2} that 
\begin{align}
&\tilde{z}_t+\lambda_1\tilde{z}_x=F(x,t)+K\lambda_1,
\quad\tilde{w}_t+\lambda_2\tilde{w}_x=F(x,t)-K\lambda_2.
\label{eqn:force3}
%\\&-M-\int^{\infty}_0b(y)dy\leqq \tilde{z},\; \tilde{w}\leqq M.\nonumber
\end{align}

\begin{figure}[htbp]
	\begin{center}
		\vspace{-2ex}
		\hspace{-6ex}
		\includegraphics[scale=0.35]{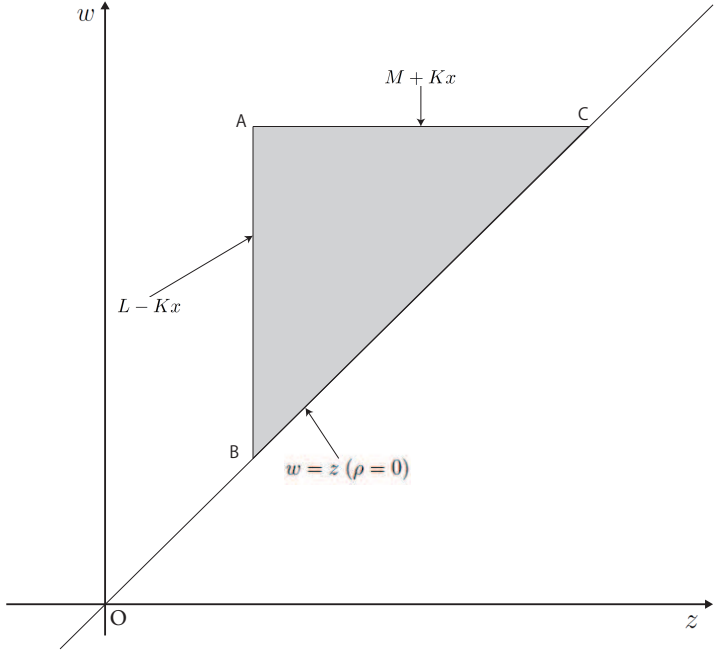}
	\end{center}
	\caption{The invariant region in $(z,w)$-plane}
	\label{Fig:estimate}
\end{figure}

Let us investigate the effects of the source term of  $\eqref{eqn:force3}$ on sides AB and AC (see Fig. \ref{Fig:estimate}), where sides AB and AC represent  $z=L-Kx$ and $w=M+Kx$ respectively. 
In these regions, $z$ and $w$ satisfy the following.
\begin{align}
\begin{split}
L-Kx\leqq z,\;\; w\leqq M+Kx.
\end{split}
\label{eqn:outline1}
\end{align}

First, we consider the source term of  $\eqref{eqn:force3}_1$ on 
the side AB. Since the following 
\begin{align*}
L-Kx=z,\;\; w\leqq M+Kx
\end{align*}
hold on the side AB, we find that 
\begin{align}
\lambda_1=v-\rho^{\theta}\geqq v-\frac{\rho^{\theta}}{\theta}=z\geqq L-Kx.
\label{eqn:outline2}
\end{align}

Therefore,    
we obtain 
\begin{align*}
\tilde{z}_t+\lambda_1\tilde{z}_x=&F(x,t)+K\lambda_1\\
\geqq&K(\lambda_1-1)\quad(\text{from $\eqref{eqn:condition-X}$})\\
\geqq&K\left(L-Kx-1\right)\hspace*{3.5ex}(\text{from $\eqref{eqn:outline2}$})\\
\geqq&0\quad(\text{from $\eqref{eqn:condition-M}$}).
\end{align*}

Second, we consider the source term of  $\eqref{eqn:force3}_2$ on 
the side AC. Since the following 
\begin{align*}
L-Kx\leqq z,\;\; w= M+Kx
\end{align*}
hold on the side AC, we find that 
\begin{align}
\lambda_2\geqq \frac{1+\theta}2M+\frac{1-\theta}2L+\theta Kx.
\label{eqn:outline3}
\end{align}

Therefore,    
we obtain 
\begin{align*}
\tilde{w}_t+\lambda_2\tilde{w}_x=&F(x,t)-K\lambda_2\\
\leqq&-K(\lambda_2-1)\quad(\text{from $\eqref{eqn:condition-X}$})\\
\leqq&-K\left( \frac{1+\theta}2M+\frac{1-\theta}2L+\theta Kx-1\right)\quad(\text{from $\eqref{eqn:outline3}$})\\
\leqq&0\quad(\text{from $\eqref{eqn:condition-M}$}).
\end{align*}

We thus conclude that the source term of $\eqref{eqn:force3}_1$ is positive on the side AB and the source term of $\eqref{eqn:force3}_2$ is negative on the side AC. We apply the maximum principle to $\tilde{z}$ and $\tilde{w}$. Then, if initial data are contained in the triangle ABC, 
we find that $\tilde{z}\geqq L$ and $\tilde{w}\leqq M$.
It follows from  that the solution remains in the same triangle. 
This implies that the following region 
\begin{align*}
\Delta_{x}=
\left\{
(z,w);\;
\rho\geqq0,\;L-Kx\leqq z,
%\right.
%\\&
%\left.
w\leqq M+Kx
\right\}
\end{align*}
is an invariant region for the initial boundary value problem $(\ref{eqn:IP})$. 
This idea 
is also used in \cite{T2}--\cite{T6}.

Although the above argument is formal, it is essential. In fact, we shall implicitly use the property of the source terms in Section 4. However, we cannot justify 
the above argument by the standard difference scheme such as Godunov or Lax-Friedrichs 
scheme (cf. \cite{DC1} and \cite{DC2}). Therefore, we introduce the modified Lax Friedrichs scheme in Section 3.

The present paper is organized as follows.
In Section 2, we review the Riemann problem and the properties of Riemann 
solutions. In Section 3, we construct approximate solutions by 
a modified Lax Friedrichs scheme. The approximate solutions consist of the steady state solutions and Riemann solutions stated in the previous section. In Section 4, we drive the bounded estimate of our approximate solutions. This section is the main point of the present paper. In Section 5, we apply 
a fixed point theorem and prove the existence of a fixed point.

\section{Preliminary}
In this section, we first review some results of the Riemann solutions 
for the homogeneous system of gas dynamics. Consider the homogeneous system 
\begin{equation}\left\{\begin{array}{ll}
\rho_t+m_x=0,\\
\displaystyle{m_t+\left(\frac{m^2}{\rho}+p(\rho)\right)_x=0,
	\quad{p}(\rho)=\rho^{\gamma}/\gamma.}
\end{array}\right.
\label{eqn:homogeneous}
\end{equation}

A pair of functions $(\eta,q):{\bf R}^2\rightarrow{\bf R}^2$ is called an 
entropy--entropy flux pair if it satisfies an identity
\begin{equation}
\nabla{q}=\nabla\eta\nabla{f}.
\label{eqn:eta-q}
\end{equation}
Furthermore, if, for any fixed ${m}/{\rho}\in(-\infty,\infty)$, $\eta$ 
vanishes on the vacuum $\rho=0$, then $\eta$ is called a {\it weak entropy}. 
For example, the mechanical energy--energy flux pair 
\begin{equation}
\eta_*:=\frac12\frac{m^2}{\rho}+\frac1{\gamma(\gamma-1)}\rho^{\gamma},\quad
q_*:=m\left(\frac12\frac{m^2}{\rho^2}+\frac{\rho^{\gamma-1}}{\gamma-1}\right)  
\label{eqn:mechanical}
\end{equation}
should be a strictly convex weak entropy--entropy flux pair. 

The jump discontinuity in a weak solutions to (\ref{eqn:homogeneous}) must satisfy the following Rankine--Hugoniot condition 
\begin{align}
\lambda(u-u_0)=f(u)-f(u_0),
\label{eqn:R-H}
\end{align}
where $\lambda$ is the propagation speed of the discontinuity, 
$u_0=(\rho_0,m_0)$ and $u=(\rho,m)$ are the corresponding left and 
right state, respectively. 
A jump discontinuity is called a {\it shock} if it satisfies the entropy 
condition 
\begin{align}
\lambda(\eta(u)-\eta(u_0))-(q(u)-q(u_0))\geqq0
\label{eqn:entropy-condition}
\end{align}
for any convex entropy pair $(\eta,q)$.

There are two distinct types of rarefaction and shock curves in the 
isentropic gases.
Given a left state $(\rho_0,m_0)$ or $(\rho_0,v_0)$, the possible states 
$(\rho,m)$ or $(\rho,v)$ that can be connected to $(\rho_0,m_0)$ or 
$(\rho_0,v_0)$ on the right by a rarefaction or a shock curve form
a 1-rarefaction wave curve $R_1(u_0)$, a 2-rarefaction wave curve $R_2(u_0)$,
a 1-shock curve $S_1(u_0)$ and a 2-shock curve $S_2(u_0)$: 
\begin{align*}
&R_1(u_0):w=w_0,\;\rho<\rho_0,\quad R_2(u_0):z=z_0,\;\rho>\rho_0,\\
&S_1(u_0):
\displaystyle{v-v_0=-\sqrt{\frac1{\rho\rho_0}\frac{p(\rho)-p(\rho_0)}
		{\rho-\rho_0}}(\rho-\rho_0)\quad\rho>\rho_0>0,}\\
&S_2(u_0):
\displaystyle{v-v_0=\sqrt{\frac1{\rho\rho_0}\frac{p(\rho)-p(\rho_0)}
		{\rho-\rho_0}}(\rho-\rho_0)\quad\rho<\rho_0,}
\end{align*}
respectively. Here we notice that shock wave curves are deduced from the\linebreak 
Rankine--Hugoniot condition (\ref{eqn:R-H}).

\begin{figure}[htbp]
	\begin{center}
		\vspace{-2ex}
		\hspace{-6ex}
		\includegraphics[scale=0.42]{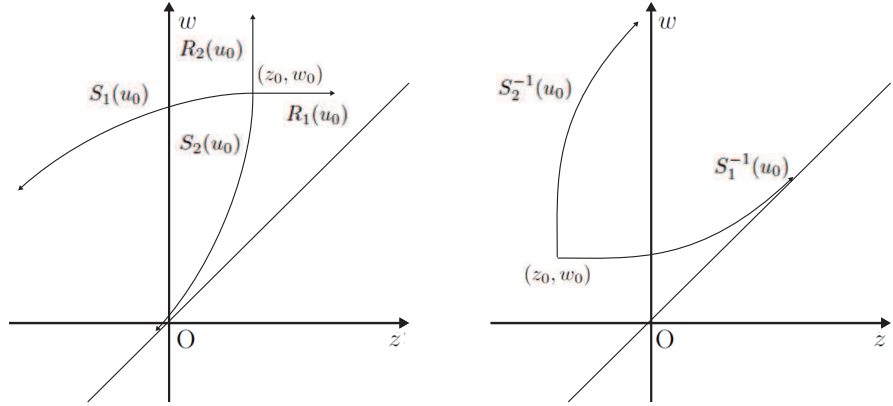}
	\end{center}
	\caption{The rarefaction curves, the shock curves and the inverse rarefaction curves in $(z,w)$-plane}
	%\label{Fig:Laval}
\end{figure}

\begin{remark}\label{rem:S-Rw}
	\normalfont
	Assume that there exists $C>1$ such that 
	\begin{align*}
	1/C\leqq\rho/\rho_0\leqq{C}.
	\end{align*}
	Then, considering $w$ along $S_1(u_0)$, we have 
	\begin{align*}
	w|_{S_1(v_0)}&=v_0-\sqrt{\frac1{\rho\rho_0}\frac{p(\rho)-p(\rho_0)}
		{\rho-\rho_0}}(\rho-\rho_0)+\frac{\rho^{\theta}}{\theta}\\
	&=w(v_0)
	+O(1)(\rho_0)^{\frac{\gamma-7}{2}}(\rho-\rho_0)^3,
	\end{align*} 
	where $O(1)$ depends only on $C$.

	Considering $z$ along $S_2(u_0)$, we similarly have
	\begin{align*}
	z|_{S_2(v_0)}&=v_0+\sqrt{\frac1{\rho\rho_0}\frac{p(\rho)-p(\rho_0)}
		{\rho-\rho_0}}(\rho-\rho_0)-\frac{\rho^{\theta}}{\theta}\\
	&=z(v_0)+O(1)
	(\rho_0)^{\frac{\gamma-7}{2}}(\rho-\rho_0)^3,
	\end{align*} 
	where $O(1)$ depends only on $C$.
	These representation show that $S_1$ (resp. $S_2$) and $R_1$ (resp. $R_2$) have a tangency of 
	second order at the point $(\rho_0,u_0)$.
\end{remark}

\subsection{Riemann Solution}
Given a right state $(\rho_0,m_0)$ or $(\rho_0,v_0)$, the possible states 
$(\rho,m)$ or $(\rho,v)$ that can be connected to $(\rho_0,m_0)$ or 
$(\rho_0,v_0)$ on the left by a shock curve constitute
1-inverse shock curve $S_1^{-1}(u_0)$ and 2-inverse shock curve 
$S_2^{-1}(u_0)$: 
\begin{align*}&&S_1^{-1}(u_0):
\displaystyle{v-v_0=-\sqrt{\frac1{\rho\rho_0}\frac{p(\rho)-p(\rho_0)}
		{\rho-\rho_0}}(\rho-\rho_0),\quad\underline{\rho<\rho_0},}\\
&&S_2^{-1}(u_0):
\displaystyle{v-v_0=\sqrt{\frac1{\rho\rho_0}\frac{p(\rho)-p(\rho_0)}
		{\rho-\rho_0}}(\rho-\rho_0),\quad\underline{\rho>\rho_0>0},}
\end{align*}
respectively.

Next we define a rarefaction shock. Given 
$u_0,u$ on $S_i^{-1}(u_0)\;(i=1,2)$, 
we call the piecewise 
constant solution to (\ref{eqn:homogeneous}), which 
consists of the left and right states $u_0,u$ a {\it rarefaction shock}. 
Here, notice the following: although the inverse shock curve has the same 
form as the shock curve, the underline expression in $S_i^{-1}(u_0)$ 
is different from the corresponding part in $S_i(u_0)$.
Therefore the rarefaction shock does not satisfy the entropy condition.

We shall use a rarefaction shock in approximating a rarefaction wave.
In particular, when we consider a rarefaction shock, we call the inverse shock 
curve connecting $u_0$ and $u$ a {\it rarefaction shock curve}.

From the properties of these curves in phase plane $(z,w)$, we can construct 
a unique solution for the Riemann problem 
\begin{equation}
u|_{t=0}=\left\{\begin{array}{ll}
u_-,\quad{x}<x_0,\\
u_+,\quad{x}>x_0,
\label{eqn:Riemann}
\end{array}\right.\end{equation}
where $x_0\in(-\infty,\infty)$, $\rho_{\pm}\geqq0$ and $m_{\pm}$ are constants 
satisfying $|m_{\pm}|\leqq{C}\rho_{\pm}$. We denote the solution 
the Riemann solution $(u_-,u_+)$.

For the Riemann problem, the following invariant region exists. 
\begin{lemma}\label{lem:invariant-region}
	For ${\rm B}_+\geqq{\rm B}_-$, the region $\sum({\rm B}_+,{\rm B}_-)
	=\{(\rho,\rho{v})\in{\bf R}^2:w=v+\rho^{\theta}/\theta,\:
	z=v-\rho^{\theta}/\theta,\:w\leqq{\rm B}_+,\: z\geqq{\rm B}_-,\:
	w-z\geqq0\}$ 
	is invariant with respect to both of the Riemann problem 
	$(\ref{eqn:Riemann})$
	and the average of the Riemann solutions in $x$. More precisely, if 
	the Riemann data lie in $\sum({\rm B}_+,{\rm B}_-)$, the corresponding Riemann 
	solutions $(\rho(x,t),m(x,t))=(\rho(x,t),\rho(x,t)v(x,t))$ 
	lie in $\sum({\rm B}_+,{\rm B}_-)$, and their corresponding 
	averages in $x$ are also in $\sum({\rm B}_+,{\rm B}_-)$, namely              
	\begin{align*}
	\left(\frac1{b-a}\int_a^b\rho(x,t)dx,\frac1{b-a}\int_a^bm(x,t)dx\right)
	\in\textstyle{\sum\nolimits({\rm B}_+,{\rm B}_-)}.
	\end{align*}
\end{lemma}
Lemma \ref{lem:invariant-region} can be found in \cite[Lemma 3.3]{C3}.

From the properties of these curves in phase plane $(z,w)$, we can construct 
a unique solution for the Riemann problem \eqref{eqn:homogeneous} and \eqref{eqn:Riemann}
\begin{equation}
u|_{t=0}=\left\{\begin{array}{ll}
u_-,\quad{x}<x_0,\\
u_+,\quad{x}>x_0,
\label{eqn:Riemann}
\end{array}\right.\end{equation}
and the Riemann initial boundary problem \eqref{eqn:homogeneous} and \eqref{eqn:RiemannIB1}
\begin{equation}
u|_{t=0}=u_+,\quad{u}|_{x=0}=u_b,\quad{x}>0,\;t>0,
\label{eqn:RiemannIB1}
\end{equation} and \eqref{eqn:homogeneous} and \eqref{eqn:RiemannIB2} 
\begin{equation}
u|_{t=0}=u_-,\quad{x}<1,\;t>0,
\label{eqn:RiemannIB2}
\end{equation} 
where $x_0\in(0,1)$, $\rho_{\pm}\geqq0$ and $m_{\pm}$ are constants 
satisfying $|m_{\pm}|\leqq{C}\rho_{\pm}$, $\lambda_1(u_{\pm})=v_{\pm}-(\rho_{\pm})^{\theta}\geqq0$ and $\lambda_1(u_b)=v_b-(\rho_b)^{\theta}\geqq0$. 

For the problem (\ref{eqn:homogeneous}) and (\ref{eqn:Riemann}), 
we can consult \cite[Subsection 3.2]{C3}.

Then the following theorem holds \cite[Theorem 3.2]{C3}. 
\begin{theorem}
	There exists a unique piecewise smooth entropy solution 
	$(\rho(x,t),$ $m(x,t))$ 
	containing the vacuum state $(\rho=0)$ on the upper plane $t>0$ for each 
	problem of {\rm (\ref{eqn:Riemann})}, {\rm (\ref{eqn:RiemannIB1})} and 
    {\rm (\ref{eqn:RiemannIB2})}
	satisfying 
	
	{\rm (1)} For the Riemann problem {\rm (\eqref{eqn:homogeneous})} and {\rm (\ref{eqn:Riemann})},
	\begin{align*}\left\{\begin{array}{lll}
	w(\rho(x,t),m(x,t))\leqq\max(w(\rho_-,m_-),w(\rho_+,m_+)),\\
	z(\rho(x,t),m(x,t))\geqq\min(z(\rho_-,m_-),z(\rho_+,m_+)),\\
	w(\rho(x,t),m(x,t))-z(\rho(x,t),m(x,t))\geqq0.
	\end{array}\right.
	\end{align*}
	
	{\rm (2)} For the Riemann initial boundary problem {\rm (\eqref{eqn:homogeneous})} and {\rm (\ref{eqn:RiemannIB1})},\begin{align*}\left\{\begin{array}{lll}
	w(\rho(x,t),m(x,t))\leqq\max(w(\rho_b,m_b),w(\rho_+,m_+)),\\
	z(\rho(x,t),m(x,t))\geqq\min(z(\rho_b,m_b),z(\rho_+,m_+)),\\
	w(\rho(x,t),m(x,t))-z(\rho(x,t),m(x,t))\geqq0.
	\end{array}\right.
	\end{align*}

    {\rm (3)} For the Riemann initial boundary problem {\rm (\ref{eqn:RiemannIB2})},
	the solution is $u_-$.
\end{theorem}

Such solutions also have the following properties:  
\begin{lemma}\label{lem:invariant-region}
	For ${\rm B}_+\geqq{\rm B}_-$, the region $\sum({\rm B}_+,{\rm B}_-)
	=\{(\rho,\rho{u})\in{\bf R}^2:w=u+\rho^{\theta}/\theta,\:
	z=u-\rho^{\theta}/\theta,\:w\leqq{\rm B}_+,\: z\geqq{\rm B}_-,\:
	w-z\geqq0\}$ 
	is invariant with respect to both of the Riemann problem {\rm (\ref{eqn:Riemann})}, the Riemann initial boundary value problem 
    \eqref{eqn:RiemannIB1} 
	and the average of the Riemann solutions in $x$. More precisely, if 
	the Riemann data lie in $\sum({\rm B}_+,{\rm B}_-)$, the corresponding Riemann 
	solutions $(\rho(x,t),m(x,t))=(\rho(x,t),\rho(x,t)u(x,t))$ 
	lie in $\sum({\rm B}_+,{\rm B}_-)$, and their corresponding 
	averages in $x$ also in $\sum({\rm B}_+,{\rm B}_-)$, namely              
	\begin{align*}
	\left(\frac1{b-a}\int_a^b\rho(x,t)dx,\frac1{b-a}\int_a^bm(x,t)dx\right)
	\in\textstyle{\sum\nolimits({\rm B}_+,{\rm B}_-)}.
	\end{align*} 
\end{lemma}
The proof of Lemma 2.2 can be found in \cite[Lemma 3.3]{C3}.

\section{Construction of Approximate Solutions}
\label{sec:construction-approximate-solutions}
In this section, we construct approximate solutions. In the strip 
$0\leqq{t}\leqq{1}$, we denote these 
approximate solutions by $u^{\varDelta}(x,t)
=(\rho^{\varDelta}(x,t),m^{\varDelta}(x,t))$. 
For $N_x\in{\bf N}$, we define the space 
lengths by ${\varDelta}x=1/(2N_x)$. 
Using $M$ in \eqref{eqn:condition-M}, we take time mesh length ${\varDelta}{x}$ such that 
\begin{align*}
\frac{{\varDelta}x}{{\varDelta}{t}}=[\hspace{-1.2pt}[2(M+K)]\hspace{-1pt}]+1,
\end{align*}
where $[\hspace{-1.2pt}[x]\hspace{-1pt}]$ is the greatest integer 
not greater than $x$. Then we define $N_t=1/(2{\varDelta t})\in{\bf N}$.
In addition, 
we set 
\begin{align*}
(j,n)\in{\bf N}_x\times{\bf N}_t,
\end{align*}
where ${\bf N}_x=\{0,1,2,\ldots,2N_x\}$ and ${\bf N}_t=\{0,1,2,\ldots,2N_t\}$.

First we define $u^{\varDelta}(x,-0)$ by 
\begin{align*}
u^{\varDelta}(x,-0)=u_0(x)
\end{align*}
and set 
\begin{align*}
J_n=\{k\in{\bf N}_x;k+n=\text{odd}\}.
\end{align*}
Then, for $j\in J_0$, we define $E_j^0(u)$ by
\begin{align*}
E_j^0(u)=\frac1{2{\varDelta}x}\int_{{(j-1)}{\varDelta}x}^{(j+1){\varDelta}x}
u^{\varDelta}(x,-0)dx.
\end{align*}

Next, assume that $u^{\varDelta}(x,t)$ is defined for $t<n{\varDelta}{t}$. 

\vspace*{1ex}

(i) $n$ is even

\vspace*{1ex}

Then, for $j\in J_n$, we define $E^n_j(u)$ by 
\begin{align*}
E^n_j(u)=\frac1{2{\varDelta}x}\int_{{(j-1)}{\varDelta}x}^{(j+1){\varDelta}x}u^{\varDelta}(x,n{\varDelta}{t}-0)dx.
\end{align*}

\vspace*{1ex}

(ii) $n$ is odd

\vspace*{1ex}

Then, for $j\in J_n\setminus\{0,2N_x\}$, we define $E^n_j(u)$ by 
\begin{align*}
E^n_j(u)=\frac1{2{\varDelta}x}\int_{{(j-1)}{\varDelta}x}^{(j+1){\varDelta}x}u^{\varDelta}(x,n{\varDelta}{t}-0)dx;
\end{align*}
for $j\in \{0,2N_x\}$, we define $E^n_j(u)$ by 
\begin{align*}
E^n_0(u)=&\frac1{{\varDelta}x}\int_{0}^{{\varDelta}x}u^{\varDelta}(x,n{\varDelta}{t}-0)dx,\\
E^n_{2N_x}(u)=&\frac1{{\varDelta}x}\int_{(2N_x-1){\varDelta}x}^{2N_x{\varDelta}x}u^{\varDelta}(x,n{\varDelta}{t}-0)dx,
\end{align*}
Moreover, for $j\in J_n$, we define $u_j^n=(\rho_j^n,m_j^n)$ in the similar manner to 
(i).

Then, for $j\in J_n$, we define $u_j^n=(\rho_j^n,m_j^n)$ as follows.\\
We choose $\delta$ such that $1<\delta<1/(2\theta)$. If 
\begin{align*}
E^n_j(\rho):=
\frac1{2{\varDelta}x}\int_{{(j-1)}{\varDelta}x}^{(j+1){\varDelta}x}\rho^{\varDelta}(x,n{\varDelta}{t}-0)dx<({\varDelta}x)^{\delta},
\end{align*} 
we define $u_j^n$ by $u_j^n=(0,0)$;
otherwise, setting\begin{align}
\begin{split}
{z}_j^n:&=
\max\left\{z(E_j^n(u)),\;L-K({j{\varDelta}x})\right\}
\\&\hspace*{100pt}\mbox{and}\\
w_j^n:&=\min\left\{w(E_j^n(u)),\;M+K({j{\varDelta}x})\right\}
,
\end{split}
\label{eqn:def-u^n_j}
\end{align}

we define $u_j^n$ by
\begin{align*}
u_j^n:=(\rho_j^n,m_j^n)
:=\left(\left\{\frac{\theta(w_j^n-z_j^n)}{2}\right\}
^{1/\theta},
\left\{\frac{\theta(w_j^n-z_j^n)}{2}\right\}^{1/\theta}
\frac{w_j^n+z_j^n}{2}\right).
\end{align*}

\begin{remark}\normalfont
	We find 
	\begin{align}
	\begin{split}
	L-K({j{\varDelta}x})\leqq z(u_j^n),\quad
	{w}(u_j^n)\leqq M+K({j{\varDelta}x}).
	\end{split}
	\label{eqn:remark3.1}
	\end{align}

	This implies that we cut off the parts where 
	$z(E_j^n(u))<L-K({j{\varDelta}x})$
	and $w(E_j^n(u))>M+K({j{\varDelta}x})$
	in  defining $z(u_j^n)$ and 
	${w}(u_j^n)$. Observing \eqref{eqn:goal}, the order of these cut parts is $o({\varDelta}x)$. The order is so small that we can deduce the compactness and convergence of our approximate solutions.
\end{remark}

\subsection{Construction of Approximate Solutions in the Cell}
\label{subsec:construction-approximate-solutions}
By using $u_j^n$ defined above, we 
construct the approximate solutions  in the cell $n{\varDelta}{t}\leqq{t}<(n+1){\varDelta}{t}\quad 
(n\in{\bf N}_t),\quad (j-1){\varDelta}x\leqq{x}<(j+1){\varDelta}x\quad
(j\in J_n\setminus\{0,1,2N_x-1,2N_x\})$. 

We first solve a Riemann problem with initial data $(u_{j-1}^n,u_{j+1}^n)$. 
Call constants $u_{\rm L}(=u_{j-1}^n), u_{\rm M}, u_{\rm R}(=u_{j+1}^n)$ the left, middle and 
right states, respectively. Then the following four cases occur.
\begin{itemize}
	\item {\bf Case 1} A 1-rarefaction wave and a 2-shock arise. 
	\item {\bf Case 2} A 1-shock and a 2-rarefaction wave arise. 
	\item {\bf Case 3} A 1-rarefaction wave and a 2-rarefaction arise.
	\item {\bf Case 4} A 1-shock and a 2-shock arise.
\end{itemize}
We then construct approximate solutions $u^{\varDelta}(x,t)$ by perturbing 
the above Riemann solutions. We consider only the case in which $u_{\rm M}$ is away from the vacuum. The other case (i.e., the case where $u_{\rm M}$ is near the vacuum) is a little technical. Therefore, we postpone the case near the vacuum to Appendix A. In addition, we omit the $L^{\infty}$ estimates for the case in this paper. The detail can be found in \cite{T2}.  

%\newpage
\vspace*{10pt}
{\bf The case where $u_{\rm M}$ is away from the vacuum}

Let $\alpha$ be a constant satisfying $1/2<\alpha<1$. Then we can choose 
a positive value $\beta$ small enough such that $\beta<\alpha$, $1/2+\beta/2<\alpha<
1-2\beta$, $\beta<2/(\gamma+5)$ and $(9-3\gamma)\beta/2<\alpha$.

We first consider the case where $\rho_{\rm M}>({\varDelta}x)^{\beta}$, 
which  means $u_{\rm M}$ is away from the vacuum. In this step, we 
consider Case 1 in particular. The constructions of Cases 2--4 are similar 
to that of Case 1.

Consider the case where a 1-rarefaction wave and a 2-shock arise as a Riemann 
solution with initial data $(u_j^n,u_{j+1}^n)$. Assume that 
$u_{\rm L},u_{\rm M}$ 
and $u_{\rm M},u_{\rm R}$ are connected by a 1-rarefaction and a 2-shock 
curve, respectively. \vspace*{10pt}\\
{\it Step 1}.\\
In order to approximate a 1-rarefaction wave by a piecewise 
constant {\it rarefaction fan}, we introduce the integer  
\begin{align*}
p:=\max\left\{[\hspace{-1.2pt}[(z_{\rm M}-z_{\rm L})/({\varDelta}x)^{\alpha}]
\hspace{-1pt}]+1,2\right\},
\end{align*}
where $z_{\rm L}=z(u_{\rm L}),z_{\rm M}=z(u_{\rm M})$ and $[\hspace{-1.2pt}[x]\hspace{-1pt}]$ is the greatest integer 
not greater than $x$. Notice that
\begin{align}
p=O(({\varDelta}x)^{-\alpha}).
\label{eqn:order-p}
\end{align}
Define \begin{align*}
z_1^*:=z_{\rm L},\;z_p^*:=z_{\rm M},\;w_i^*:=w_{\rm L}\;(i=1,\ldots,p),
\end{align*}
and 
\begin{align*}
z_i^*:=z_{\rm L}+(i-1)({\varDelta}x)^{\alpha}\;(i=1,\ldots,p-1).
%\label{eqn:z^*_i^*-z^*_{i+1}}
\end{align*}
We next introduce the rays $x=(j+1/2){\varDelta}x+\lambda_1(z_i^*,z_{i+1}^*,w_{\rm L})
(t-n{\varDelta}{t})$ separating finite constant states 
$(z_i^*,w_i^*)\;(i=1,\ldots,p)$, 
where  
\begin{align*}
\lambda_1(z_i^*,z_{i+1}^*,w_{\rm L}):=v(z_i^*,w_{\rm L})
-S(\rho(z_{i+1}^*,w_{\rm L}),\rho(z_i^*,w_{\rm L})),
\end{align*}
\begin{align*}
\rho_i^*:=\rho(z_i^*,w_{\rm L}):=\left(\frac{\theta(w_{\rm L}-z_i^*)}2\right)^{1/\theta}\;,
\quad{v}_i^*:={v}(z_i^*,w_{\rm L}):=\frac{w_{\rm L}+z_i^*}2
\end{align*}
and

\begin{align}
S(\rho,\rho_0):=\left\{\begin{array}{lll}
\sqrt{\displaystyle{\frac{\rho(p(\rho)-p(\rho_0))}{\rho_0(\rho-\rho_0)}}}
,\quad\mbox{if}\;\rho\ne\rho_0,\\
\sqrt{p'(\rho_0)},\quad\mbox{if}\;\rho=\rho_0.
\end{array}\right.
\label{eqn:s(,)}
\end{align}

We call this approximated 1-rarefaction wave a {\it 1-rarefaction fan}.

\vspace*{10pt}
{\it Step 2}.\\
In this step, we replace the above constant states  with the following functions of $x$:

\begin{definition}\label{def:steady-state}\normalfont

	For given constants $x_d$ satisfying  $(j-1){\varDelta}x\leqq{x_d}\leqq(j+1){\varDelta}x$ and 
	\begin{align}
	\begin{split}
	(z_d,w_d):=\left(\frac{m_d}{\rho_d}-\frac{(\rho_d)^{\theta}}
	{\theta},\frac{m_d}{\rho_d}+\frac{(\rho_d)^{\theta}}{\theta}\right)
	\quad\text{or}\quad
	u_d=(\rho_d,m_d)
	\end{split}
	\label{eqn:steady-state-data}
	\end{align}
	satisfying $|m_d|\leqq{C}\rho_d$, we set \begin{align*}
	z(x)=z_d-K(x-{x_d}),\quad
	w(x)=w_d+K(x-{x_d}).
	\end{align*}

	Using $z(x)$ and $w(x)$, we define  
	\begin{align}
	u(x)=(\rho(x),m(x))
	\label{eqn:solution-steady-state}
	\end{align}
	by the relation 
	\eqref{eqn:relation-Riemann} as follows.
	\begin{align}
	v(x)=\frac{w(x)+z(x)}2,\;\rho(x)=\left(\frac{\theta(w(x)-z(x))}2\right)^{1/\theta},\;m(x)=\rho(x) v(x).
	\label{eqn:relation-Riemann2}
	\end{align}
	We then define $\bar{\mathcal U}(x,x_d,u_d)$ with data $u_d$ at $x_d$ as  
	\eqref{eqn:solution-steady-state} .
\end{definition}

Moreover, for given functions $\bar{u}(x)$, we define $z(x,t)$ and $w(x,t)$ by \begin{align}
\begin{split}
z^{\varDelta}(x,t)=&\bar{z}^{\varDelta}(x)
+\left\{F(x,t)+K\lambda_1(\bar{u}^{\varDelta}(x))\right\}
(t-n{\varDelta{t}}),\\
w^{\varDelta}(x,t)=&\bar{w}^{\varDelta}(x)+
\left\{F(x,t)-K\lambda_2(\bar{u}^{\varDelta}(x))\right\}(t-n{\varDelta{t}}).
\end{split}
\label{eqn:fractional-step}
\end{align}

Then, using $z(x,t)$ and $w(x,t)$, we define $u(x,t)=(\rho(x,t),m(x,t))$ in a similar manner to
\eqref{eqn:relation-Riemann2}.  We denote $u(x,t)$ by ${\mathcal U}(x,t;\bar{u})$.

Let $\bar{u}_{\rm L}(x)$ and $\bar{u}_{\rm R}(x)$ be $\bar{\mathcal U}
(x,(j-1){\varDelta}x,u_{\rm L})$ and 
$\bar{\mathcal U}(x,(j+1){\varDelta}x,u_{\rm R})$, respectively. Set 
$
\bar{u}_1(x):=\bar{u}_{\rm L}(x),\;u_1(x,t)={\mathcal U}(x,t;\bar{u}_1),\;u_{\rm R}(x,t)={\mathcal U}(x,t;\bar{u}_{\rm R})$ and $x_1:=(j-1){\varDelta}x.
$

\begin{figure}[htbp]
	\begin{center}
		\hspace{-2ex}
		\includegraphics[scale=0.3]{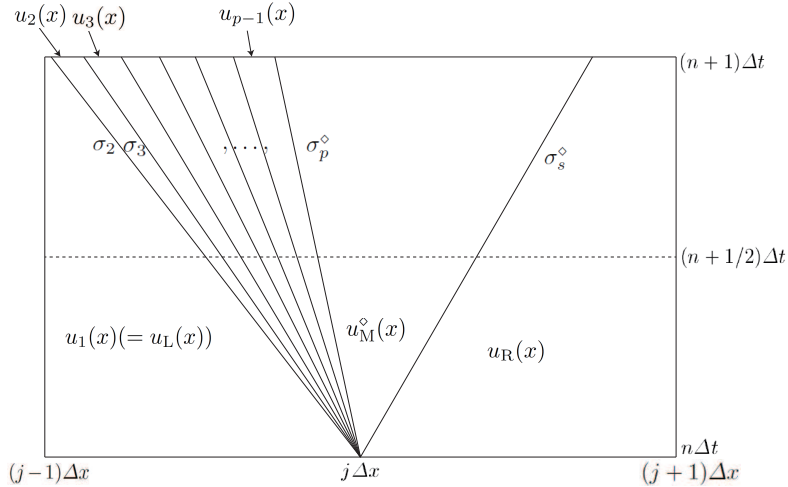}
	\end{center}
	\caption{The approximate solution in the case where a 1-rarefaction and 
		a 2-shock arise in the cell.}
	\label{case1-1cell}
\end{figure}

First, by the implicit function theorem, we determine a propagation speed $\sigma_2$ and $u_2=(\rho_2,m_2)$ such that 
\begin{itemize}
	\item[(1.a)] $z_2:=z(u_2)=z^*_2$
	\item[(1.b)] the speed $\sigma_2$, the left state ${u}_1(x_2,(n+1/2){\varDelta}t)$ and the right state $u_2$ satisfy the Rankine--Hugoniot conditions, i.e.,
	\begin{align*}
	f(u_2)-f({u}_1(x_2,(n+1/2){\varDelta}t))=\sigma_2(u_2-{u}_1(x_2,(n+1/2){\varDelta}t)),
	\end{align*}
\end{itemize}
where $x_2:=j{\varDelta}x+\sigma_2{\varDelta}t/2$. Then we fill up by ${u}_1(x)$ the sector where $n{\varDelta}t\leqq{t}<
(n+1){\varDelta}t,(j-1){\varDelta}x\leqq{x}<{j}{\varDelta}x+
\sigma_2(t-n{\varDelta}t)$ (see Figure \ref{case1-1cell}) 
and set $\bar{u}_2(x)=\bar{\mathcal U}(x,x_2,u_2)$ and $u_2(x,t)={\mathcal U}(x,t;\bar{u}_2)$.

Assume that $u_k$, ${u}_k(x,t)$ and a propagation speed $\sigma_k$ with
$\sigma_{k-1}<\sigma_k$ are defined. Then we similarly determine
$\sigma_{k+1}$ and $u_{k+1}=(\rho_{k+1},m_{k+1})$ such that 
\begin{itemize}
	\item[($k$.a)] $z_{k+1}:=z(u_{k+1})=z^*_{k+1}$,
	\item[($k$.b)] $\sigma_{k}<\sigma_{k+1}$,
	\item[($k$.c)] the speed 
	$\sigma_{k+1}$, 
	the left state ${u}_k(x_{k+1},(n+1/2){\varDelta}t)$ and the right state $u_{k+1}$ satisfy 
	the Rankine--Hugoniot conditions, 
\end{itemize}
where 
$x_{k+1}:=j{\varDelta}x+\sigma_{k+1}{\varDelta}t/2$. Then we fill up by ${u}_k(x,t)$ the sector where
$n{\varDelta}t\leqq{t}<(n+1){\varDelta}t,{j}{\varDelta}x+\sigma_k(t-{\varDelta}t)\leqq{x}<{j}{\varDelta}x+\sigma_{k+1}(t-n{\varDelta}t)$  (see Figure \ref{case1-1cell}) and 
set $\bar{u}_{k+1}(x)=\bar{\mathcal U}(x,x_{k+1},u_{k+1})$ and $u_{k+1}(x,t)={\mathcal U}(x,t;\bar{u}_{k+1})$.  
By induction, we define $u_i$, ${u}_i(x,t)$ and $\sigma_i$ $(i=1,\ldots,p-1)$.
Finally, we determine a propagation speed $\sigma_p$ and $u_p=(\rho_p,m_p)$ such that
\begin{itemize}
	\item[($p$.a)] $z_p:=z(u_p)=z^*_p$,
	\item[($p$.b)] the speed $\sigma_p$, 
	and the left state ${u}_{p-1}(x_p,(n+1/2){\varDelta}t)$ and the right state $u_p$ satisfy the Rankine--Hugoniot conditions, 
\end{itemize}where $x_p:=j{\varDelta}x+\sigma_p{\varDelta}t/2$. 
We then fill up by ${u}_{p-1}(x,t)$ and $u_p$ the sector where
$n{\varDelta}t\leqq{t}<(n+1){\varDelta}t,{j}{\varDelta}x+\sigma_{p-1}
(t-n{\varDelta}t)\leqq{x}<{j}{\varDelta}x+\sigma_{p}(t-n{\varDelta}t)$ 
and the line $n{\varDelta}t\leqq{t}<(n+1){\varDelta}t,x={j}{\varDelta}x+\sigma_{p}(t-n{\varDelta}t)$, respectively.

Given $u_{\rm L}$ and $z_{\rm M}$ with $z_{\rm L}\leqq{z}_{\rm M}$, we denote 
this piecewise functions of $x$ 1-rarefaction wave by 
$R_1^{\varDelta}(z_{\rm M})(u_{\rm L})$. Notice that from the construction 
$R^{\varDelta}_1(z_{\rm M})(u_{\rm L})$ connects $u_{\rm L}$ and $u_p$ 
with $z_p=z_{\rm M}$. 

Now we fix ${u}_{\rm R}(x,t)$ and ${u}_{p-1}(x,t)$. 
Let $\sigma_s$ be 
the propagation speed of the 2-shock connecting $u_{\rm M}$ and $u_{\rm R}$.
Choosing ${\sigma}^{\diamond}_p$ near to $\sigma_p$, ${\sigma}^{\diamond}_s$ 
near to 
$\sigma_s$ and $u^{\diamond}_{\rm M}$ near to $u_{\rm M}$, we fill up by ${u}^{\diamond}_{\rm M}(x,t)=
{\mathcal U}(x,t,\bar{u}^{\diamond}_{\rm M})$ the gap between $x=j{\varDelta}x+{\sigma}^{\diamond}_{p}
(t-n{\varDelta}{t})$ and $x=j{\varDelta}x+{\sigma}^{\diamond}_s(t-n{\varDelta}{t})$, such that 
\begin{itemize}
	\item[(M.a)] $\sigma_{p-1}<\sigma^{\diamond}_p<\sigma^{\diamond}_s$, 
	\item[(M.b)] the speed ${\sigma}^{\diamond}_p$, the left and right states 
	${u}_{p-1}(x^{\diamond}_{p},(n+1/2){\varDelta}t),{u}^{\diamond}_{\rm M}(x^{\diamond}_{p},(n+1/2){\varDelta}t)$ 
	satisfy the Rankine--Hugoniot conditions,
	\item[(M.c)] the speed ${\sigma}^{\diamond}_s$, the left and right 
	states ${u}^{\diamond}_{\rm M}(x^{\diamond}_{s},(n+1/2){\varDelta}t),{u}_{\rm R}(x^{\diamond}_{s},(n+1/2){\varDelta}t)$ satisfy the Rankine--Hugoniot conditions, 
\end{itemize}
where $\bar{u}^{\diamond}_{\rm M}(x)=
\bar{\mathcal U}(x,j{\varDelta}x,u^{\diamond}_{\rm M})$,\;$x^{\diamond}_{p}:=j{\varDelta}x+\sigma^{\diamond}_{p}{\varDelta}
/2$ and $x^{\diamond}_s:=j{\varDelta}x+\sigma^{\diamond}_s{\varDelta}
/2$. 

We denote this approximate Riemann solution, which consists of \eqref{eqn:solution-steady-state}, by ${u}^{\varDelta}(x,t)$. The validity of the above construction is demonstrated in \cite[Appendix A]{T2}.

\begin{remark}\label{rem:middle-time}\normalfont
	${u}^{\varDelta}(x,t)$ satisfies the Rankine--Hugoniot conditions
	at the middle time of the cell, $t_{\rm M}:=(n+1/2){\varDelta}t$.
\end{remark}

%\vspace*{10pt}

\begin{remark}\label{rem:approximate}\normalfont
	The approximate solution $u^{\varDelta}(x,t)$ is piecewise smooth in each of the 
	divided parts of the cell. Then, in the divided part, $u^{\varDelta}(x,t)$ satisfies
	\begin{align*}
	(u^{\varDelta})_t+f(u^{\varDelta})_x-g(x,u^{\varDelta})=O(\varDelta x).
	\end{align*}
\end{remark}

\subsection{Construction of Approximate Solutions near the Boundary}
For $n\in{\bf N}_t$, separating four parts, we construct approximate solutions near the boundary as follows.
\begin{enumerate}
	\item Near the left boundary $x=0$ 
	\begin{enumerate}
		\item $0\leqq{x}<2{\varDelta}x,\; n{\varDelta}{t}\leqq{t}<(n+1){\varDelta}{t}$, $n$ is even

		We first solve a Riemann problem with initial data 
		\begin{align*}
		u|_{t=n{\varDelta}t}=
		\begin{cases}
		u_b,\quad{x}<0,\\
		u^n_1,\quad{x}>0,
		\end{cases}\end{align*}For the Riemann solution, we 
		construct the approximate solutions in the cell $-2{\varDelta}x\leqq{x}<2{\varDelta}x,\; n{\varDelta}{t}\leqq{t}<(n+1){\varDelta}{t}$ in the same manner to Section 3.1. Then we 
		define the resultant approximate solution in $0\leqq{x}<2{\varDelta}x,\; n{\varDelta}{t}\leqq{t}<(n+1){\varDelta}{t}$ by $u^{\varDelta}(x,t)$ in this area.	
		\item $0\leqq{x}<{\varDelta}x,\; n{\varDelta}{t}\leqq{t}<(n+1){\varDelta}{t}$, $n$ is odd

			We first solve a Riemann problem with initial data 
		\begin{align*}
		u|_{t=n{\varDelta}t}=
		\begin{cases}
		u_b,\quad{x}<0,\\
		u^n_0,\quad{x}>0,
		\end{cases}\end{align*}For the Riemann solution, we 
		construct the approximate solutions in the cell $-{\varDelta}x\leqq{x}<{\varDelta}x,\; n{\varDelta}{t}\leqq{t}<(n+1){\varDelta}{t}$ in the same manner to Section 3.1. Then we 
		define 
		the resultant approximate solution in $0\leqq{x}<{\varDelta}x,\; n{\varDelta}{t}\leqq{t}<(n+1){\varDelta}{t}$ by $u^{\varDelta}(x,t)$ in this area.
	\end{enumerate}
\item Near the right boundary $x=1$ 
\begin{enumerate}
\item $1-2{\varDelta}x\leqq{x}<1,\; n{\varDelta}{t}\leqq{t}<(n+1){\varDelta}{t}$, $n$ is even

We first solve a Riemann problem with initial data 
\begin{align*}
u|_{t=n{\varDelta}t}=
\begin{cases}
u^n_1,\quad{x}<1,\\
u^n_1,\quad{x}>1,
\end{cases}\end{align*}For the Riemann solution, we 
construct the approximate solutions in the cell $1-2{\varDelta}x\leqq{x}<1+2{\varDelta}x,\; n{\varDelta}{t}\leqq{t}<(n+1){\varDelta}{t}$ in the same manner to Section 3.1. Then we 
define the resultant approximate solution in $1-2{\varDelta}x\leqq{x}<1,\; n{\varDelta}{t}\leqq{t}<(n+1){\varDelta}{t}$ by $u^{\varDelta}(x,t)$ in this area.	
\item $1-{\varDelta}x\leqq{x}<1,\; n{\varDelta}{t}\leqq{t}<(n+1){\varDelta}{t}$, $n$ is odd

We first solve a Riemann problem with initial data 
\begin{align*}
u|_{t=n{\varDelta}t}=
\begin{cases}
u^n_0,\quad{x}<1,\\
u^n_0,\quad{x}>1,
\end{cases}\end{align*}For the Riemann solution, we 
construct the approximate solutions in the cell $1-{\varDelta}x\leqq{x}<1+{\varDelta}x,\; n{\varDelta}{t}\leqq{t}<(n+1){\varDelta}{t}$ in the same manner to Section 3.1. Then we 
define the resultant approximate solution in $1-{\varDelta}x\leqq{x}<1,\; n{\varDelta}{t}\leqq{t}<(n+1){\varDelta}{t}$ by $u^{\varDelta}(x,t)$ in this area.	
\end{enumerate}
\end{enumerate}

\section{$L^{\infty}$ Estimate of the Approximate Solutions}\label{sec:bound}
We estimate Riemann invariants of $u^{\varDelta}(x,t)$ to use the invariant 
region theory. Our aim in this section is to deduce from (\ref{eqn:remark3.1}) the following
theorem:\begin{theorem}\label{thm:bound}
	\begin{align}
	\begin{split}
	\displaystyle L-Kx-{\it o}({\varDelta}x)
	\leqq {z}^{\varDelta}(x,(n+1)\varDelta t-0),\\
	\displaystyle {w}^{\varDelta}(x,(n+1)\varDelta t-0)
	\leqq M+Kx+{\it o}({\varDelta}x),
	\end{split}
	\label{eqn:goal}
	\end{align}
	where ${\it o}({\varDelta}x)$ depends only on $M$.
\end{theorem}

In this section, we first assume \begin{align}
M\geqq\left(1+K\right)+\varepsilon
\label{eqn:condition-M2}
\end{align}	
instead of \eqref{eqn:condition-M}, where $\varepsilon$ is any fixed positive 
value.

\vspace*{5pt}
Throughout this paper, by the Landau symbols such as $O({\varDelta}x)$,
$O(({\varDelta}x)^2)$ and $o({\varDelta}x)$,
we denote quantities whose moduli satisfy a uniform bound depending only on
$M$ and $X$ unless we specify otherwise.  In addition, for simplicity, we 
denote $w(\bar{u}_i(x))$ and $z(\bar{u}_i(x))$ by $\bar{w}_i(x)$ and 
$\bar{z}_i(x)$.

Now, in the previous section, we have constructed 
$u^{\varDelta}(x,t)$. Then, the following four cases occur.
\begin{itemize}
	\item In {\bf Case 1}, the main difficulty is to obtain $(\ref{eqn:goal})_1$
	along $R^{\varDelta}_1$.  
	\item In {\bf Case 2}, the main difficulty is to obtain $(\ref{eqn:goal})_2$
	along $R^{\varDelta}_2$. 
	\item In {\bf Case 3}, (\ref{eqn:goal}) follows that of {Case 1}
	and Case 2. 
	\item In {\bf Case 4}, (\ref{eqn:goal}) is easier than that of the other 
	cases. 
\end{itemize}

Thus we treat  Case 1 in particular. In Case 1, we derive $(\ref{eqn:goal})_1$ along $R^{\varDelta}_1$   
and estimate the other quasi-steady state solutions.
We can estimate the other cases in a fashion similar to Case 1.

\subsection{Estimates of $\bar{u}^{\varDelta}(x,t)$ in Case 1}
In this step, we estimate ${u}^{\varDelta}(x,t)$ in Subsection \ref{subsec:construction-approximate-solutions}.  In this case, each component of $\bar{u}^{\varDelta}(x,t)$ 
has the following properties, which is  proved in \cite[Appendix A]{T2}:
\begin{align}
&\bullet\;\sigma_i<\sigma_{i+1}\;(i=1,\ldots,p-2),
\sigma_{p-1}<\sigma^{\diamond}_p
<\sigma^{\diamond}_s.\\
&\bullet\;\rho_i>({\varDelta}x)^{\beta}/2\;(i=1,\ldots,p-1).
\label{eqn:lower-rho}\\
%\\
&\bullet\;\mbox{Given data $z_i:=z(u_i)$ and $w_i:=w(u_i)$ at $x=x_i$, $\bar{u}_i(x,t)
	=\bar{\mathcal U}(x,x_i,u_i)$,}\nonumber\\
&\mbox{ $(i=1,\ldots,p-1)$ that is,}\nonumber\\
&(\bar{z}_{i}(x),\;\bar{w}_{i}(x))
%\\&
=\left(z_i-K(x-{x_i}),\;w_i+K(x-{x_i})\right)
\label{eqn:bar-v-Delta}\\
&\bullet\;
\bar{w}_{i+1}(x_{i+1})=w_{i+1}=\bar{w}_i(x_{i+1})
+{\it O}(({\varDelta}x)^{3\alpha-(\gamma-1)\beta})
%\nonumber\\&\hspace{30ex}
\quad
(i=1,\ldots,p-2).
%&\hspace{50ex}
\label{eqn:w_i-w_{i+1}}\\
&\bullet\;|u_{\rm M}^{\diamond}-u_{\rm M}|=O(({\varDelta}x)^{1-\frac{\gamma+1}{2}\beta}).
\label{eqn:v_m}\\
&\bullet\;\bar{u}^{\diamond}_{\rm M}(x)=\bar{\mathcal U}(x,(j+1/2){\varDelta}x,u^{\diamond}_{\rm M}).\nonumber
\\
&\bullet\;\mbox{${u}_i(x_{i+1})$ and ${u}_{i+1}(x_{i+1},(n+1/2){\varDelta}t)$
	are connected}\nonumber
\\
&\text{by a 1-rarefaction shock curve}\;\;(i=1,\ldots,p-2).\nonumber
\\
&\bullet\;\bar{w}^{\diamond}_{\rm M}(x^{\diamond}_{p})=
\bar{w}_{p-1}(x^{\diamond}_p)
+{\it O}(({\varDelta}x)^{3\alpha+(\gamma-7)\beta/2}).
\label{eqn:w_m-w_{p-1}}
\\
%\end{align}
%\begin{align}
&\bullet\;\mbox{${u}_{p-1}(x^{\diamond}_p,(n+1/2){\varDelta}t)$ and ${u}^{\diamond}_{\rm M}(x^{\diamond}_p,(n+1/2){\varDelta}t)$
	are connected by }\nonumber\\
&\;\;\mbox{a 1-shock or a 1-rarefaction shock curve.}\nonumber\\
&\bullet\;\mbox{${u}^{\diamond}_{\rm M}(x^{\diamond}_s,(n+1/2){\varDelta}t)$ and ${u}_{\rm R}(x^{\diamond}_s,(n+1/2){\varDelta}t)$
	are connected by }\nonumber\\
&\;\;\mbox{a 2-shock or a 2-rarefaction shock curve.}\nonumber
\end{align}

Now we derive (\ref{eqn:goal}) in the interior cell
$n{\varDelta}{t}\leqq{t}<(n+1){\varDelta}{t},\;(j-1){\varDelta}x\leqq{x}
\leqq(j+1){\varDelta}x$. To do this, we first consider components of 
$\bar{u}^{\varDelta}(x,t)$.

\vspace*{5pt}\hspace*{-20pt}
{\bf Estimate of $\bar{z}_i(x)\;(i=1,\ldots,p-1)$}.\\

Recalling that $z_1=z_{\rm L}=z_{j-1}^n\geqq L-K({(j-1){\varDelta}x})$, 
we have 
\begin{align*}
\bar{z}_1(x)=z_1-K({x}-{(j-1){\varDelta}x})
\geqq L-K{x}.
\end{align*}

On the other hand, 
the construction of $\bar{u}_i(x)$ implies that $z_i=z_1+(i-1)({\varDelta}x)^{\alpha}$. If $i\geqq2$, since $\alpha<1$, it follows that 
\begin{align*}
z_i&
=z_1+(i-1)({\varDelta}x)^{\alpha}
\geqq L-K({(j-1){\varDelta}x})+({\varDelta}x)^{\alpha}
%\nonumber\\&
\geqq L-K({(j-1){\varDelta}x}).
%\label{eqn:estimate-z_i3}
\end{align*} 
We thus obtain 
\begin{align}
\bar{z}_i(x)=z_i-K({x}-{x_i})\geqq L-Kx.
\label{eqn:estimate-z_i}
\end{align}

{\bf Estimate of $\bar{z}_{\rm R}(x)$}.\\
Recall that $z_{\rm R}=z_{j+1}^n\geqq L-K((j+1){\varDelta}x)$. We then have 
\begin{align}
\bar{z}_{\rm R}(x)=z_{\rm R}-K({x}-{(j+1){\varDelta}x})
\geqq L-Kx.
\label{eqn:estimate-z_r}
\end{align}

{\bf Estimate of $\bar{z}^{\diamond}_{\rm M}(x)$}.\\
If ${u}^{\diamond}_{\rm M}(x^{\diamond}_s,(n+1/2){\varDelta}t)$ and ${u}_{\rm R}(x^{\diamond}_s,(n+1/2){\varDelta}t)$
are connected by a 2-shock curve, we find 
\begin{align*}
{z}^{\diamond}_{\rm M}(x^{\diamond}_s,(n+1/2){\varDelta}t)\geqq{z}_{\rm R}(x^{\diamond}_s,(n+1/2){\varDelta}t).
%\label{eqn:estimate-z_m1}
\end{align*}
Then, in view of \eqref{eqn:fractional-step}, we obtain  
\begin{align*}
\bar{z}^{\diamond}_{\rm M}(x^{\diamond}_s)\geqq\bar{z}_{\rm R}(x^{\diamond}_s).
%\label{eqn:estimate-z_m1}
\end{align*}
Therefore, from (\ref{eqn:estimate-z_r}), we have
\begin{align}
\bar{z}^{\diamond}_{\rm M}(x^{\diamond}_s)
\geqq L-K{x^{\diamond}_s}.
\label{eqn:estimate-z_m1}
\end{align}

On the other hand, we consider the case where ${u}^{\diamond}_{\rm M}(x^{\diamond}_s,(n+1/2){\varDelta}t)$ and ${u}_{\rm R}(x^{\diamond}_s,(n+1/2){\varDelta}t)$ are connected by a 2-rarefaction shock curve. 
First, recall that $u_{\rm M}$ and $u_{\rm R}$ are connected, not by a 
2-rarefaction shock curve but by a 2-shock curve. Since $|{u}^{\diamond}_{\rm M}(x^{\diamond}_s,(n+1/2){\varDelta}t)-\bar{u}^{\diamond}_{\rm M}(x^{\diamond}_s)|=O({\varDelta}x)$, $|\bar{u}^{\diamond}_{\rm M}(x^{\diamond}_s)-
u^{\diamond}_{\rm M}|=O({\varDelta}x)$, $|{u}_{\rm R}(x^{\diamond}_s,(n+1/2){\varDelta}t)-\bar{u}_{\rm R}(x^{\diamond}_s)|=O({\varDelta}x)$ and $|\bar{u}_{\rm R}(x^{\diamond}_s)-u_{\rm R}|=O({\varDelta}x)$, we then deduce from (\ref{eqn:v_m}) that 
\begin{align*}
|\bar{u}^{\diamond}_{\rm M}(x^{\diamond}_s)-\bar{u}_{\rm R}(x^{\diamond}_s)|=
O(({\varDelta}x)^{1-(\gamma+1)\beta/2}).
\end{align*}
Therefore, from Remark \ref{rem:S-Rw} and the fact that $\beta<2/(\gamma+5)$, we conclude that  
\begin{align}
\bar{z}^{\diamond}_{\rm M}(x^{\diamond}_s)&=\bar{z}_{\rm R}
(x^{\diamond}_s)-O(({\varDelta}x)
^{3(1-\frac{\gamma+1}2\beta)+\frac{\gamma-7}2\beta})
\geqq L-K{x^{\diamond}_s}-o({\varDelta}x).
\label{eqn:estimate-z_m2}
\end{align}
Therefore, from \eqref{eqn:estimate-z_m1}--\eqref{eqn:estimate-z_m2}, we obtain 
\begin{align}
\bar{z}^{\diamond}_{\rm M}(x)
&=\bar{z}^{\diamond}_{\rm M}(\bar{x}^{\diamond}_s)+K({x}-{x^{\diamond}_s})
\geqq L-Kx-o({\varDelta}x).
\label{eqn:estimate-z_m3}
\end{align}

\vspace*{5pt}\hspace*{-20pt}{\bf Estimate of $\bar{w}_i(x)\;(i=1,\ldots,p-1)$}.
\\

First, we recall that 
\begin{align*}
w_1=w(u_{j-1}^n)\leqq M+K((j-1){\varDelta}x).
\end{align*}
We then assume that 
\begin{align}
w_i\leqq M+K{x_i}+i\cdot O(({\varDelta}x)^{3\alpha-(\gamma-1)\beta}).
\label{eqn:estimate-w_i1}
\end{align}
It follows that  
\begin{align*}
\bar{w}_i(x)=w_i+K(x-{x_i})\leqq M+K{x}
+i\cdot O(({\varDelta}x)^{3\alpha-(\gamma-1)\beta}).
\end{align*}
From \eqref{eqn:w_i-w_{i+1}}, we obtain 
\begin{align}
{w}_{i+1}\leqq M+K{x_{i+1}}+(i+1)\cdot O(({\varDelta}x)^{3\alpha-(\gamma-1)\beta}).
\label{eqn:estimate-w_i2}
\end{align}
Therefore, \eqref{eqn:estimate-w_i1} holds for any $i$.

In view of  \eqref{eqn:order-p} and \eqref{eqn:estimate-w_i1}, since $3\alpha-(\gamma-1)\beta>1$, we thus conclude that  
\begin{align}
\bar{w}_i(x)=w_i+K({x}-{x_i})\leqq M+K{x}+o({\varDelta}x).
\label{eqn:estimate-w_i3}
\end{align}

\vspace*{5pt}\hspace*{-20pt}{\bf Estimate of $\bar{w}^{\diamond}_{\rm M}(x)$}.
\\
Combining the fact that $(9-3\gamma)\beta/2<\alpha$, (\ref{eqn:w_m-w_{p-1}}) and (\ref{eqn:estimate-w_i3}), we thus have 
\begin{align}
\bar{w}^{\diamond}_{\rm M}(x)&=\bar{w}^{\diamond}_{\rm M}(x^{\diamond}_p)
=\bar{w}_{p-1}(x^{\diamond}_p)+{\it o}({\varDelta}x)\leqq M+K{x}+o({\varDelta}x). 
\label{eqn:estimate-w_m}
\end{align}

\vspace*{5pt}\hspace*{-20pt}{\bf Estimate of $\bar{w}_{\rm R}(x)$}.
\\
Recalling $w_{\rm R}=w(u^n_{j+1})\leqq
M+K((j+1){\varDelta}x)$, it follows
that 
\begin{align}
\bar{w}_{\rm R}(x)=w_{\rm R}+K({x}-{(j+1){\varDelta}x})
\leqq M+K{x}.
\label{eqn:estimate-w_r}
\end{align}

%\newpage
%\vspace*{10pt}\hspace*{-20pt}
{\bf Estimate of $u^{\varDelta}(x,t)$ in the interior cell $n{\varDelta}{t}\leqq{t}<(n+1){\varDelta}{t},\;(j-1){\varDelta}x\leqq{x}<(j+1){\varDelta}x$}.

Let us derive $\eqref{eqn:goal}_1$.

\begin{figure}[htbp]
	\begin{center}
		\vspace{-2ex}
		\hspace{-6ex}
		\includegraphics[scale=0.4]{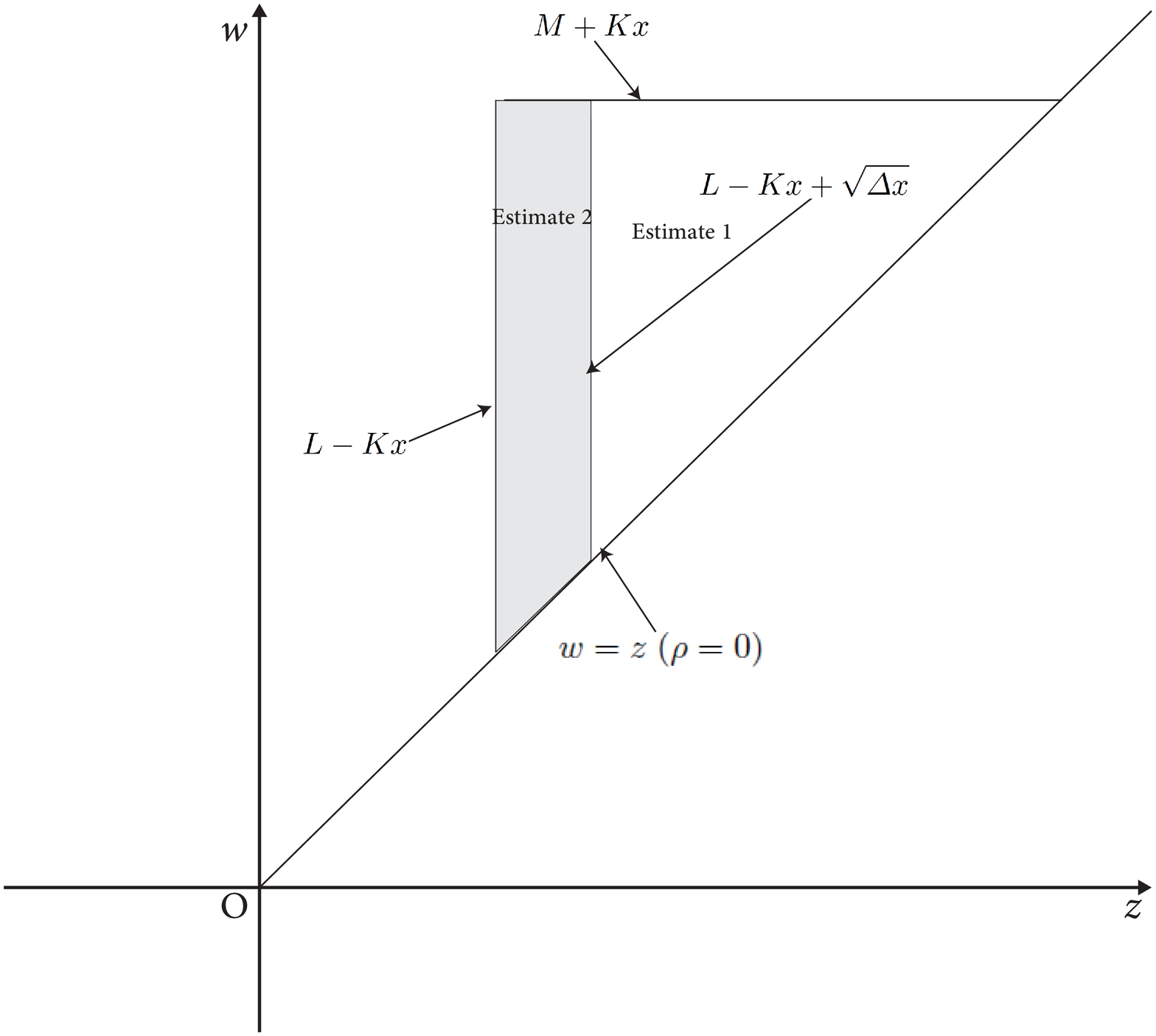}
	\end{center}
	\caption{The invariant region in $(z,w)$-plane}
	%\label{Fig:Laval}
\end{figure}
{\bf Estimate 1}

We first consider the case where $z(\bar{u}^{\varDelta}(x,t))\geqq L-K{x}+
\sqrt{{\varDelta}x}$. 
On the other hand, we recall that  $z({u}^{\varDelta}(x,t))=z(\bar{u}^{\varDelta}(x,t))+{\it O}({\varDelta}x)$. Therefore, choosing ${\varDelta}x$ small enough, we conclude $\eqref{eqn:goal}_1$.

{\bf Estimate 2}

We next consider the case where 
\begin{align}
z(\bar{u}^{\varDelta}(x,t))<L-K{x}+\sqrt{{\varDelta}x}.
\label{eqn:estimate-z2}
\end{align} 
This case is the validity of Subsection 1.1. Recalling its argument, let us deduce $\eqref{eqn:goal}_1$.

From \eqref{eqn:estimate-z_i}, 
\eqref{eqn:estimate-z_r}, \eqref{eqn:estimate-z_m3}, \eqref{eqn:estimate-w_i3}, 
\eqref{eqn:estimate-w_m} and \eqref{eqn:estimate-w_r}, we have 
\begin{align}
\begin{split}
z(\bar{u}^{\varDelta}(x,t))\geqq L-Kx-{\it o}({\varDelta}x),\quad
w(\bar{u}^{\varDelta}(x,t))\leqq M+K{x}+{\it o}({\varDelta}x).
\end{split}
\label{eqn:estimate-z3}
\end{align}

Then, from \eqref{eqn:estimate-z2}--\eqref{eqn:estimate-z3}, choosing ${\varDelta}x$ small enough, we find that 
\begin{align}
\lambda_1(\bar{u}^{\varDelta}(x))\geqq L-K{x}+O(\sqrt{{\varDelta}x}).
\label{eqn:lambda_1}
\end{align}

Then, from \eqref{eqn:condition-X} and \eqref{eqn:condition-M2}, we obtain 
\begin{align*}
z({u}^{\varDelta}&(x,t))\nonumber\\
&=z(\bar{u}^{\varDelta}(x,t))
+\left\{F(x,t)
+K\lambda_1(\bar{u}^{\varDelta}(x))\right\}
(t-n{\varDelta{t}})\nonumber\\
&\geqq z(\bar{u}^{\varDelta}(x,t))
+K\left(\lambda_1(\bar{u}^{\varDelta}(x))-1\right)
(t-n{\varDelta{t}})\nonumber\\
&\geqq z(\bar{u}^{\varDelta}(x,t))+K\left\{L-K{x}-1+O(\sqrt{{\varDelta}x})\right\}
(t-n{\varDelta{t}})\nonumber\\
&\geqq z(\bar{u}^{\varDelta}(x,t)).
%\label{eqn:estimate-z5}
\end{align*}
As a result, 
from \eqref{eqn:estimate-z3}, 
we drive $\eqref{eqn:goal}_1$. We can similarly obtain $\eqref{eqn:goal}_2$.

\section{Recurrence formula}From Remark \ref{rem:approximate}, $u^{\varDelta}$ satisfy 
\begin{align*}
(u^{\varDelta})_t+f(u^{\varDelta})_x-g(x,t,u^{\varDelta})=O(\varDelta x)
\end{align*}
on the divided part in the cell where $u^{\varDelta}$ are smooth. Moreover, $u^{\varDelta}$ satisfy an entropy condition (see \cite[Lemma 5.1--Lemma 5.4]{T2}) along 
discontinuous lines approximately. Then, applying the Green formula to $(u^{\varDelta})_t+f(u^{\varDelta})_x-g(x,t,u^{\varDelta})$ in the cell $(j-1){\varDelta}x\leqq{x}<(j+1){\varDelta}x,\;n{\varDelta}{t}\leqq{t}<(n+1){\varDelta}{t}\quad(j\in J_{n+1},\;n\in{\bf N}_t)$, we have
\begin{align}
\begin{split}
\rho^{n+1}_j=&\frac{\rho^n_{j+1}+\rho^n_{j-1}}2-\frac{\varDelta t}{2\varDelta x}\left\{m^n_{j+1}-m^n_{j-1}\right\}\\
&+R(x_{j+1},t_n,u^n_{j+1})-R(x_{j-1},t_n,u^n_{j-1})+o(\varDelta x), 
\\
%\end{split}
%\label{eqn:density}\\
%\begin{split}
m^{n+1}_j=&\frac{m^n_{j+1}+m^n_{j-1}}2-\frac{\varDelta t}{2\varDelta x}\left\{\frac{(m^n_{j+1})^2}{\rho^n_{j+1}}+p(\rho^n_{j+1})
-\frac{(m^n_{j-1})^2}{\rho^n_{j-1}}-p(\rho^n_{j-1})\right\}\\
&+S(x_{j+1},t_n,u^n_{j+1})-S(x_{j-1},t_n,u^n_{j-1})\\
%&+\frac{1}{2\varDelta x}\int^{t_{n+1}}_{t_n}\int^{x_{j+1}}_{x_{j-1}}F(x,t)
%{\rho^{\varDelta}(x,t)}dxdt+O((\varDelta x)^2)\\
&+\frac{F(x_{j+1},t_n)\rho^n_{j+1}+F(x_{j-1},t_n)\rho^n_{j-1}}2{\varDelta t}+o(\varDelta x),
\end{split}
\label{eqn:recurrence1}
\end{align}
where $t_n=n{\varDelta t}$ and 
\begin{align*}
R(x,t,u)=&-\frac{\varDelta x}{4}K(\rho)^{1-\theta}
+\frac{(\varDelta t)^2}{4\varDelta x}
\left\{F(x,t)\rho-K\rho^{\theta+1}+K\frac{m}{\rho^{\theta+1}}\right\},\\
S(x,t,u)=&-{\varDelta t}\hspace*{0.2ex}F(x,t)m
-\frac{\varDelta x}{4}K\frac{m}{\rho^{\theta}}
\\&
+\frac{(\varDelta t)^2}{4\varDelta x}
\left\{2F(x,t)m-K\rho^{\theta}m+K\frac{m^3}{\rho^{\theta+2}}\right\}.
\end{align*}
Moreover, from \eqref{eqn:goal} and \cite[Lemma 6.1]{T3}, we have
\begin{align}
 L-K(j{\varDelta}x)-{\it o}({\varDelta}x)
\leqq {z}(u^n_j),\quad
 {w}(u^n_j)
\leqq M+K(j{\varDelta}x)+{\it o}({\varDelta}x),\quad \rho^n_j\geqq0
\label{eqn:goal2}
\end{align}

Then, we define a sequence $\tilde{u}^n_{j}=(\tilde{\rho}^{n}_j,\tilde{m}^{n}_j)$ as follows.

\begin{align}
\begin{split}
\tilde{\rho}^{n+1}_j=&\frac{\tilde{\rho}^n_{j+1}+\tilde{\rho}^n_{j-1}}2-\frac{\varDelta t}{2\varDelta x}\left\{\tilde{m}^n_{j+1}-\tilde{m}^n_{j-1}\right\}\\
&+R(x_{j+1},t_n,\tilde{u}^n_{j+1})-R(x_{j-1},t_n,\tilde{u}^n_{j-1}),
\\
%\end{split}
%\label{eqn:density2}\\
%\begin{split}
\tilde{m}^{n+1}_j=&\frac{\tilde{m}^n_{j+1}+\tilde{m}^n_{j-1}}2-\frac{\varDelta t}{2\varDelta x}\left\{\frac{(\tilde{m}^n_{j+1})^2}{\tilde{\rho}^n_{j+1}}+p(\tilde{\rho}^n_{j+1})
-\frac{(\tilde{m}^n_{j-1})^2}{\tilde{\rho}^n_{j-1}}-p(\tilde{\rho}^n_{j-1})\right\}\\
&+S(x_{j+1},t_n,\tilde{u}^n_{j+1})-S(x_{j-1},t_n,\tilde{u}^n_{j-1})\\
%&+\frac{1}{2\varDelta x}\int^{t_{n+1}}_{t_n}\int^{x_{j+1}}_{x_{j-1}}F(x,t)
%{\rho^{\varDelta}(x,t)}dxdt+O((\varDelta x)^2)\\
&+\frac{F(x_{j+1},t_n)\tilde{\rho}^n_{j+1}+F(x_{j-1},t_n)\tilde{\rho}^n_{j-1}}2{\varDelta t},\quad(j\in J_{n+1},\;n\in {\bf N}_t)
\end{split}
\label{eqn:recurrence2}\\
&\hspace*{-7ex}\tilde{u}^0_{j}={u}^0_{j} \quad (j\in J_0).
\end{align}

Therefore, from \eqref{eqn:recurrence1}--\eqref{eqn:recurrence2}, there exists $\delta({\varDelta}x)>0$ satisfying $\delta({\varDelta}x)\rightarrow0$ as 
${\varDelta}x\rightarrow 0$, such that 
\begin{align}
\begin{split}
&L-K(j{\varDelta}x)-\delta({\varDelta}x)
\leqq {z}(\tilde{u}^n_j),\quad
{w}(\tilde{u}^n_j)
\leqq M+K(j{\varDelta}x)+\delta({\varDelta}x),\quad\tilde{\rho}^n_j\geqq0.
\end{split}
\label{eqn:goal3}
\end{align}

Then, we define a map $F:{\bf R}^{4N_x+2}\rightarrow{\bf R}^{4N_x+2}$ as follows.
\begin{align}
&\left(\left\{{z}(\tilde{u}^0_{j})\right\}^{2N_x}_{j=0},
\left\{{w}(\tilde{u}^0_j)\right\}^{2N_x}_{j=0}\right)\mapsto
%\\
&\left(\left\{{z}(\tilde{u}^n_{j})+\delta({\varDelta}x)\right\}^{2N_x}_{j=0},
\left\{{w}(\tilde{u}^n_j)-\delta({\varDelta}x)\right\}^{2N_x}_{j=0}\right).
\label{eqn:map}
\end{align}
From \eqref{eqn:recurrence2}, $F$ is continuous. Moreover, from 
\eqref{eqn:goal3}, $F$ is the map from a bounded set to the same bounded set. Therefore, applying 
the Brouwer fixed point theorem to $F$, we have a fixed point $(\tilde{u}^0_j)^*$. We supply $(\tilde{u}^0_j)^*$ as initial data ${u}^0_j$ for our approximate solutions.

The following proposition and theorem can be proved in the same manner to 
\cite{T2}--\cite{T4}. 
\begin{proposition}\label{pro:compact}
	The measure sequence
	\begin{align*}
	\eta(u^{\varDelta})_t+q(u^{\varDelta})_x
	\end{align*}
	lies in a compact subset of $H_{\rm loc}^{-1}(\Omega)$ for all weak entropy 
	pair $(\eta,q)$, where $\Omega\subset{\bf R}\times{\bf R}_+$ is any bounded
	and open set. 
\end{proposition}
\begin{theorem} 
	Assume that the approximate solutions $(\rho^{\varDelta},m^{\varDelta})$ satisfy
	Theorem \ref{thm:bound} and Proposition \ref{pro:compact}. Then there is a convergent subsequence 
	in the approximate solutions $(\rho^{\varDelta}(x,t),m^{\varDelta}(x,t))$ such that
	\begin{align*}      
	&(\rho^{\varDelta_n}(x,t),m^{\varDelta_n}(x,t))\rightarrow(\rho(x,t),m(x,t))
	\hspace{2ex}
	\text{\rm a.e.} \;(x,t)\in(0,1)\times(0,1),\quad\\&(\varDelta x)_n\rightarrow0\;\text{as\;\;}n\rightarrow \infty.
	\end{align*} 
	The function $(\rho(x,t),m(x,t))$ is a time periodic entropy solution
	of the initial boundary value problem \eqref{eqn:IP}.
\end{theorem}

We have proved Theorem \ref{thm:main} under the condition  \eqref{eqn:condition-M2}.
However, since $\varepsilon$ are arbitrary, we conclude Theorem \ref{thm:bound} under the condition \eqref{eqn:condition-M}.

\section{Open problem for a time-periodic outer force}

In this section, we introduce some open problem related to \eqref{eqn:IC}. The most difficult point 
of the present problems is to prove that initial data and the corresponding solutions after one period 
are contained in the same bounded set. To solve this problem, we have introduced the generalized invariant region 
depending on the space variable. Observing \eqref{eqn:bound}, we find that the upper bound $U(x)$ and lower bound $L(x)$ are monotone functions (In \eqref{eqn:bound}, $U(x)=M+Kx$ and $L(x)=L-Kx$ in particular.). 
This method is used in \cite{T2}--\cite{T9} and these bounds are also monotone functions. 
The method is useful for the Cauchy problem or the initial boundary problem for a half space. 
Since we need not supply the boundary condition at $x=1$ for the present paper, our problem is the almost same as the initial boundary value problem for a half space. We can thus construct the invariant region.

However, unfortunately, our method is not useful for the periodic boundary condition or the initial boundary value problem for a bounded interval. Considering the periodic boundary condition, $U$ and $L$ need to be periodic. Next, considering the boundary condition $m|_{x=0}=m|_{x=1}=0$ instead of \eqref{eqn:BC}, $U(0)\geq-L(0)$ and $U(1)\leq-L(1)$ need to hold. Recalling our upper and lower bound are monotone, $U(x)$ and $L(x)$ cannot satisfy these conditions. 
On the other hand, as mentioned above, we can construct invariant regions for the Cauchy problem and the initial boundary value problem for a half space. 
However, since their regions are not finite, we can not apply the Brouwer fixed point theorem like \eqref{eqn:map}. Therefore, the existence of a
time periodic solution for these problems is still open.

\appendix

\section{Construction of Approximate Solutions near the vacuum}

In this step, we consider the case where $\rho_{\rm M}\leqq({\varDelta}x)^{\beta}$,
which means that $u_{\rm M}$ is near the vacuum. In this case, we cannot construct 
approximate solutions in a similar fashion to Subsection 3.1. Therefore, we must
define $u^{\varDelta}(x,t)$ in a different way.

In this appendix, we define our approximate solutions in the cell $(j-1){\varDelta}x\leqq{x}<(j+1){\varDelta}x,\;n{\varDelta}{t}\leqq{t}<(n+1){\varDelta}{t}\quad(j\in{\bf Z},\:n\in{\bf Z}_{\geqq0})$. 
We set $L_j:=L-K((j+1){\varDelta}x)$ and 
$U_j:={M}+K((j-1){\varDelta}x)$.

\vspace*{5pt}
{\bf Case 1} A 1-rarefaction wave and a 2-shock arise.

In this case, we notice that $\rho_{\rm R}\leqq ({\varDelta}x)^{\beta},\;
z_{\rm R}\geqq L_j$ and $w_{\rm R}\leqq U_j$.
\vspace*{5pt}

\begin{center}
	{\bf Definition of $\bar{u}^{\varDelta}$ in Case 1}
\end{center}

\vspace*{5pt}
{\bf Case 1.1}
$\rho_{\rm L}>({\varDelta}x)^{\beta}$

We denote $u^{(1)}_{\rm L}$ a state satisfying $ w(u_{\rm L}^{(1)})=w(u_{\rm L})$ and 
$\rho^{(1)}_{\rm L}=({\varDelta}x)^{\beta}$. 
Let $u^{(2)}_{\rm L}$ be a state connected to $u_{\rm L}$ on the right by 
$R_1^{\varDelta}(z^{(1)}_{\rm L})(u_{\rm L})$. We set 
\begin{align*}
(z^{(3)}_{\rm L},w^{(3)}_{\rm L})=
\begin{cases}
(z^{(2)}_{\rm L},w^{(2)}_{\rm L}),\quad\text{if $z^{(2)}_{\rm L}\geqq L_j$},\\
(L_j,w^{(2)}_{\rm L}+L_j-z^{(2)}_{\rm L}),\quad\text{if $z^{(2)}_{\rm L}< L_j$}.
\end{cases}
\end{align*}

Then, we define 
\begin{align*}
\bar{u}^{\varDelta}(x,t)
=\begin{cases}
R_1^{\varDelta}(z^{(1)}_{\rm L})(u_{\rm L}),\quad
\text{if $(j-1){\varDelta}x
	\leqq{x}\leqq{j}{\varDelta}x+\lambda_1(u^{(2)}_{\rm L})(t-n{\varDelta}t)$}\\\text{and $n{\varDelta}t\leqq{t}<(n+1){\varDelta}t$},\vspace*{2ex}\\
\text{a Riemann solution  $(u^{(3)}_{\rm L}$, $u_{\rm R})$},\quad\text{if ${j}{\varDelta}x+\lambda_1(u^{(2)}_{\rm L})(t-n{\varDelta}t)$}\\\text{$<x
	\leqq (j+1){\varDelta}x$ and $n{\varDelta}t\leqq{t}<(n+1){\varDelta}t$}.
\end{cases}
\end{align*}

%\vspace*{-1ex}
\begin{figure}[htbp]
	\begin{center}
		\vspace{-4ex}
		\hspace{2ex}
		\includegraphics[scale=0.32]{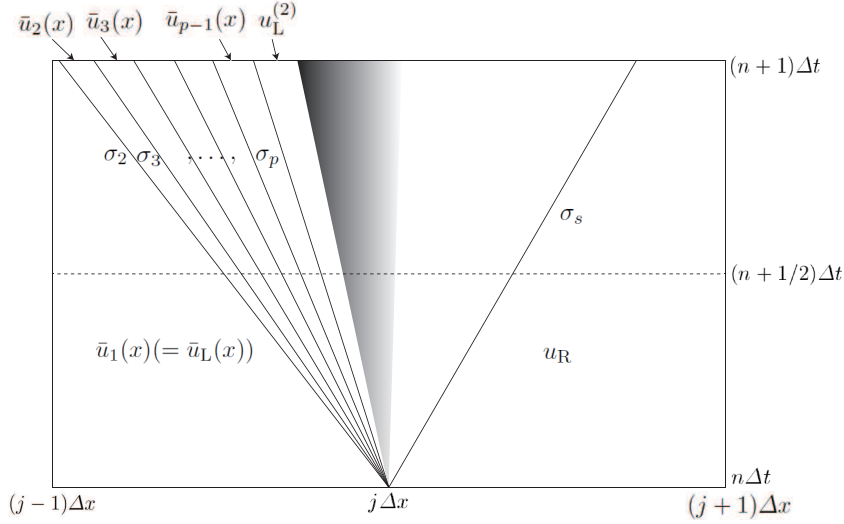}
	\end{center}
	\caption{{\bf Case 1.1}: The approximate solution $\bar{u}^{\varDelta}$ in the cell.}
	\label{Fig:case1.1(ii)}
\end{figure} 
%\vspace*{-3ex}

%\vspace*{-5ex}
{\bf Case 1.2} $\rho_{\rm L}\leqq({\varDelta}x)^{\beta}$

\vspace*{5pt}
(i) $z(u_{\rm L})\geqq{L}_j$\\
In this case, we define $u^{\varDelta}(x,t)$ as a Riemann solution 
$(u_{\rm L},u_{\rm R})$.

\vspace*{5pt}
(ii) $z(u_{\rm L})<L_j$\\
In this case, recalling $z(u_{\rm L})=z(u^n_j)\geqq L-K((j-1){\varDelta}x)$, 
we can choose $x^{(4)}$ such that $(j-1){\varDelta}x\leqq{x}^{(4)}
\leqq(j+1){\varDelta}x$ and 
%\begin{align*}
$
z(u_{\rm L})-K({x^{(4)}}-{x_{\rm L}})=L_j,
$
%\end{align*}
where $x_{\rm L}:=(j-1){\varDelta}x$.
We set   
\begin{align*}
z^{(4)}_{\rm L}:=z_{\rm L}-K({x^{(4)}}-{x_{\rm L}}),\quad{w}^{(4)}_{\rm L}:
=w_{\rm L}+K({x^{(4)}}-{x_{\rm L}}).
%\label{eqn:v^{5*}}
\end{align*}

In the region where $(j-1){\varDelta}x
\leqq{x}\leqq{j}{\varDelta}x+\lambda_1(u^{(4)}_{\rm L})(t-n{\varDelta}t)$ and 
$n{\varDelta}t\leqq{t}<(n+1){\varDelta}t$,
we define $\bar{u}^{\varDelta}(x,t)$ as a solution of (\ref{eqn:solution-steady-state}) such that 
$\bar{u}^{\varDelta}((j-1){\varDelta}x)=u_{\rm L}$.
We next solve a Riemann problem
$(u^{(4)}_{\rm L},u_{\rm R})$. In the region where ${j}{\varDelta}x+\lambda_1(u^{(4)}_{\rm L})(t-n{\varDelta}t)\leqq{x}\leqq(j+1){\varDelta}x$ and 
$n{\varDelta}t\leqq{t}<(n+1){\varDelta}t$,  
we define $\bar{u}^{\varDelta}(x,t)$ as this Riemann solution.

\vspace*{5pt}

\begin{center}
	{\bf Definition of ${u}^{\varDelta}$}
\end{center}

\vspace*{5pt}

{\bf Case 1} 
In the region 
where $\bar{u}^{\varDelta}(x,t)$ is the Riemann solution, we 
define $u^{\varDelta}(x,t)$ by $u^{\varDelta}(x,t)=\bar{u}^{\varDelta}(x,t)$; otherwise, the definition of 
$u^{\varDelta}(x,t)$ is similar to Subsection 3.1. 
Thus, for a Riemann solution near the vacuum, we define our approximate solution
as the Riemann solution itself.

\vspace*{10pt}
{\bf Case 2} A 1-shock and a 2-rarefaction wave arise.

From symmetry, this case reduces to Case 1.

\vspace*{10pt}
{\bf Case 3} A 1-rarefaction wave and a 2-rarefaction wave arise.

For $u_{\rm L}$ of Case 1, we define $u^*_{\rm L}$ and $\lambda^*_{\rm L}$ as follows. 
\begin{align*}
u^*_{\rm L}=\begin{cases}
u^{(3)}_{\rm L},\quad\text{Case 1.1},\\
u^{(4)}_{\rm L},\quad\text{Case 1.2},
\end{cases}\quad 
\lambda^*_{\rm L}=\begin{cases}
\lambda_1(u^{(2)}_{\rm L}),\quad\text{Case 1.1},\\
\lambda_1(u^{(4)}_{\rm L}),\quad\text{Case 1.2}.
\end{cases}
\end{align*}
where $\lambda_1(u)$ be the 1-characteristic speed of $u$. Then, for $u_{\rm L}$ of Case 3, we can determine $u^*_{\rm L}$ and 
$\lambda^*_{\rm L}$ in a similar manner to Case 1. 
From symmetry, for $u_{\rm R}$ of Case 3, we can also 
determine $u^*_{\rm R}$ and $\lambda^*_{\rm R}$.

In the 
region 
$(j-1){\varDelta}x\leqq{x}\leqq{j}{\varDelta}x+\lambda^*_{\rm L}(t-n{\varDelta}t),\;
{j}{\varDelta}x+\lambda^*_{\rm R}(t-n{\varDelta}t)\leqq{x}\leqq(j+1){\varDelta}x$ and 
$n{\varDelta}t\leqq{t}<(n+1){\varDelta}t$, we define $\bar{u}^{\varDelta}$ in a similar manner to Subsection 3.1. In the other 
region, we define $\bar{u}^{\varDelta}$ as the Riemann solution 
$(u^*_{\rm L},u^*_{\rm R})$.

We define ${u}^{\varDelta}$ in the same way as Case 1.

\vspace*{10pt}
{\bf Case 4} A 1-shock and a 2-shock arise.

We notice that $z_{\rm L}\geqq L_j,\;w_{\rm L}\leqq U_j,\;z_{\rm R}\geqq L_j$ and $w_{\rm R}\leqq U_j$.
In this case, we define $u^{\varDelta}(x,t)$ as the Riemann
solution $(u_{\rm L},u_{\rm R})$.
\vspace*{2ex}

We complete the construction of our approximate solutions.

\section*{Acknowledgements} 
The author would appreciate Prof. Shigeharu Takeno for his useful comments.

\end{document}